\documentclass{article}

\usepackage[left=1.25in,top=1.25in,right=1.25in,bottom=1.25in,head=1.25in]{geometry}
\usepackage{amsfonts,amsmath,amssymb,amsthm}
\usepackage{array,verbatim,float,url}
\usepackage{graphicx,subfigure,psfrag}
\usepackage{dsfont,bm}
\usepackage{natbib}

\usepackage{floatpag}
\floatpagestyle{empty}

\newtheorem{theorem}{Theorem}
\newtheorem{lemma}{Lemma}
\newtheorem{corollary}{Corollary}

\theoremstyle{remark}
\newtheorem{remark}{Remark}
\theoremstyle{definition}
\newtheorem{definition}{Definition}

\newcommand{\argmin}{\mathop{\mathrm{argmin}}}

\newcommand{\st}{\mathop{\mathrm{subject\,\,to}}}

\def\R{\mathbb{R}}
\def\E{\mathbb{E}}

\def\Cov{\mathrm{Cov}}
\def\half{\frac{1}{2}}
\def\sign{\mathrm{sign}}
\def\supp{\mathrm{supp}}
\def\rank{\mathrm{rank}}

\def\tr{\mathrm{tr}}

\def\hmu{\hat{\mu}}
\def\hbeta{\hat{\beta}}

\def\T{^T}
\def\cA{\mathcal{A}}

\def\cD{\mathcal{D}}

\def\df{\mathrm{df}}
\def\sdf{\mathrm{sdf}}
\def\hbetalas{\hbeta^\mathrm{lasso}}
\def\hbetasub{\hbeta^\mathrm{subset}}

\def\cAlas{\cA^\mathrm{lasso}}
\def\cAsub{\cA^\mathrm{subset}}
\def\hmulas{\hmu^\mathrm{lasso}}
\def\hmusub{\hmu^\mathrm{subset}}

\begin{document}

\title{Degrees of Freedom and Model Search}
\author{Ryan J. Tibshirani}
\date{}
\maketitle

\begin{abstract}
Degrees of freedom is a fundamental concept in statistical modeling,
as it provides a quantitative description of the amount of fitting
performed by a given procedure.  But, despite this fundamental
role in statistics, its behavior is not completely well-understood,
even in somewhat basic settings.  For example, it may seem
intuitively obvious that the best subset selection fit with subset
size $k$ has degrees of freedom larger than $k$, but this has not been
formally verified, nor has is been precisely studied.  At large,
the current paper is motivated by this problem, and we
derive an exact expression for the degrees of freedom of best subset
selection in a restricted setting (orthogonal predictor variables).  
Along the way, we develop a concept
that we name ``search degrees of freedom''; intuitively, for adaptive
regression procedures that perform variable selection, this is a part
of the (total) degrees of freedom that we attribute entirely to the
model selection mechanism. Finally, we establish a modest extension of
Stein's formula to cover discontinuous functions, and discuss  
its potential role in degrees of freedom and search degrees of freedom 
calculations. 
\end{abstract}

\section{Introduction}
\label{sec:intro}



Suppose that we are given observations $y\in\R^n$ from the model
\begin{equation}
\label{eq:model}
y = \mu + \epsilon, \;\;\; \text{with}\;\, 
\E(\epsilon)=0, \, \Cov(\epsilon)=\sigma^ 2 I,
\end{equation}
where $\mu\in\R^n$ is some fixed, true mean parameter of interest, and  
$\epsilon \in \R^n$ are uncorrelated errors, with zero mean and 
common marginal variance $\sigma^2 > 0$. For a function $f : \R^n
\rightarrow \R^n$,  thought of as a procedure for producing fitted
values, $\hmu = f(y)$, recall that the {\it degrees of freedom}
of $f$ is defined as \citep{bradbiased,gam}:
\begin{equation}
\label{eq:df}
\df(f) = \frac{1}{\sigma^2}\sum_{i=1}^n \Cov\big(f_i(y),y_i\big).
\end{equation}
Intuitively, the quantity $\df(f)$ reflects the effective number of
parameters used by $f$ in producing the fitted output
$\hmu$. Consider linear regression, for example, where $f(y)$ is the
least squares fit of $y$ onto predictor variables $x_1,\ldots x_p \in
\R^n$:  for this procedure $f$, our intuition gives the right
answer, as its degrees of freedom is simply $p$, the number of
estimated regression  
coefficients.\footnote{This is assuming linear independence of
  $x_1,\ldots x_p$; in general, it is the dimension of
  $\mathrm{span}\{x_1,\ldots x_p\}$.}   This, e.g., leads to
an unbiased estimate of the risk of the linear regression fit, via   
Mallows's $C_p$ criterion \citep{mallows}.

In general, characterizations of degrees of freedom are highly
relevant for purposes like model comparisons and model selection;  
see, e.g.,   \citet{bradbiased,gam,lassodf2}, and Section
\ref{sec:opt}, for more motivation.  Unfortunately, however, 
counting degrees of freedom can become quite complicated for    
nonlinear, adaptive procedures.  (By nonlinear, we mean $f$
being nonlinear as a function of $y$.)
Even for many basic adaptive procedures, explicit answers are not
known.  A good example is best subset selection, in which, for a fixed 
integer $k$, we regress on the
subset of $x_1,\ldots x_p$ of size at most $k$ giving the best
linear fit of $y$ (as measured by the residual sum of squares).
Is the degrees of freedom here larger than $k$? It seems that
the answer should be ``yes'', because even though there are $k$
coefficients in the final linear model, the variables in this model were
chosen adaptively (based on the data).  And if the answer is
indeed ``yes'', then the natural follow-up question is: how much
larger is it?  That is, how many effective parameters does it ``cost''
to search through the space of candidate models? The goal of this
paper is to investigate these questions, and related ones. 

\subsection{A motivating example}
\label{sec:motivex}

We begin by raising an interesting point: though it seems certain that
a procedure like best subset selection would suffer an inflation
of degrees of freedom, not all adaptive regression procedures do.
In particular, the {\it lasso} \citep{lasso,bp}, which also performs
variable selection in the linear model setting, presents a very
different story in terms of its degrees of freedom.  Stacking
the predictor variables $x_1,\ldots x_p$ along the columns of 
a matrix $X \in \R^{n\times p}$, the lasso estimate can be
expressed as:
\begin{equation}
\label{eq:lasso}
\hbetalas = \argmin_{\beta \in \R^p} \,
\half\|y-X\beta\|_2^2 + \lambda\|\beta\|_1,
\end{equation}
where $\lambda\geq 0$ is a tuning parameter, controlling the level
of sparsity.  Though not strictly necessary for our discussion, we
assume for simplicity that $X$ has columns in general position, which
ensures uniqueness of the lasso solution \smash{$\hbetalas$}
(see, e.g., \citet{lassounique}).
We will write $\cAlas \subseteq \{1,\ldots p\}$ to denote the indices of 
nonzero coefficients in \smash{$\hbetalas$}, called the
support or active set of \smash{$\hbetalas$}, also expressed as
\smash{$\cAlas=\supp(\hbetalas)$}. 

The lasso admits a simple formula for its degrees of freedom.
\begin{theorem}[\citealt{lassodf,lassodf2}]
\label{thm:lassodf}
Provided that the variables (columns) in $X$ are in general
position, the lasso fit \smash{$\hmu^\mathrm{lasso}=X\hbetalas$} has
degrees of freedom 
\begin{equation*}
\df(\hmu^\mathrm{lasso}) = \E|\cAlas|,   
\end{equation*}
where $|\cAlas|$ is the size of the lasso active set 
\smash{$\cAlas=\supp(\hbetalas)$}.
The above expectation assumes that $X$ and $\lambda$ are fixed, and is 
taken over the sampling distribution $y \sim N(\mu,\sigma^2 I)$.  
\end{theorem}

In other words, the degrees of freedom of the lasso fit is the number
of selected variables, in expectation. This is somewhat remarkable
because, as with subset selection, the lasso uses the data to choose
which variables to put in the model. So how can its degrees of freedom
be equal to the (average) number of selected variables, and not more?  
The key realization is that the lasso shrinks the coefficients of
these variables towards zero, instead of perfoming a full least
squares fit. This shrinkage is due to the $\ell_1$ penalty that
appears in \eqref{eq:lasso}. Amazingly, the ``surplus'' from
adaptively building the model is exactly accounted for by the
``deficit'' from shrinking the coefficients, so that altogether (in
expectation), the degrees of freedom is simply the number of variables 
in the model.    

\begin{remark} An analogous result holds for an entirely arbitrary
predictor matrix $X$ (not necessarily having columns in
general position), see \citet{lassodf2}; analogous results also exist for the
generalized lasso problem (special cases of which are the fused lasso
and trend filtering), see \citet{genlasso,lassodf2}.
\end{remark}

\begin{figure}[htb]
\centering
\includegraphics[width=0.55\textwidth]{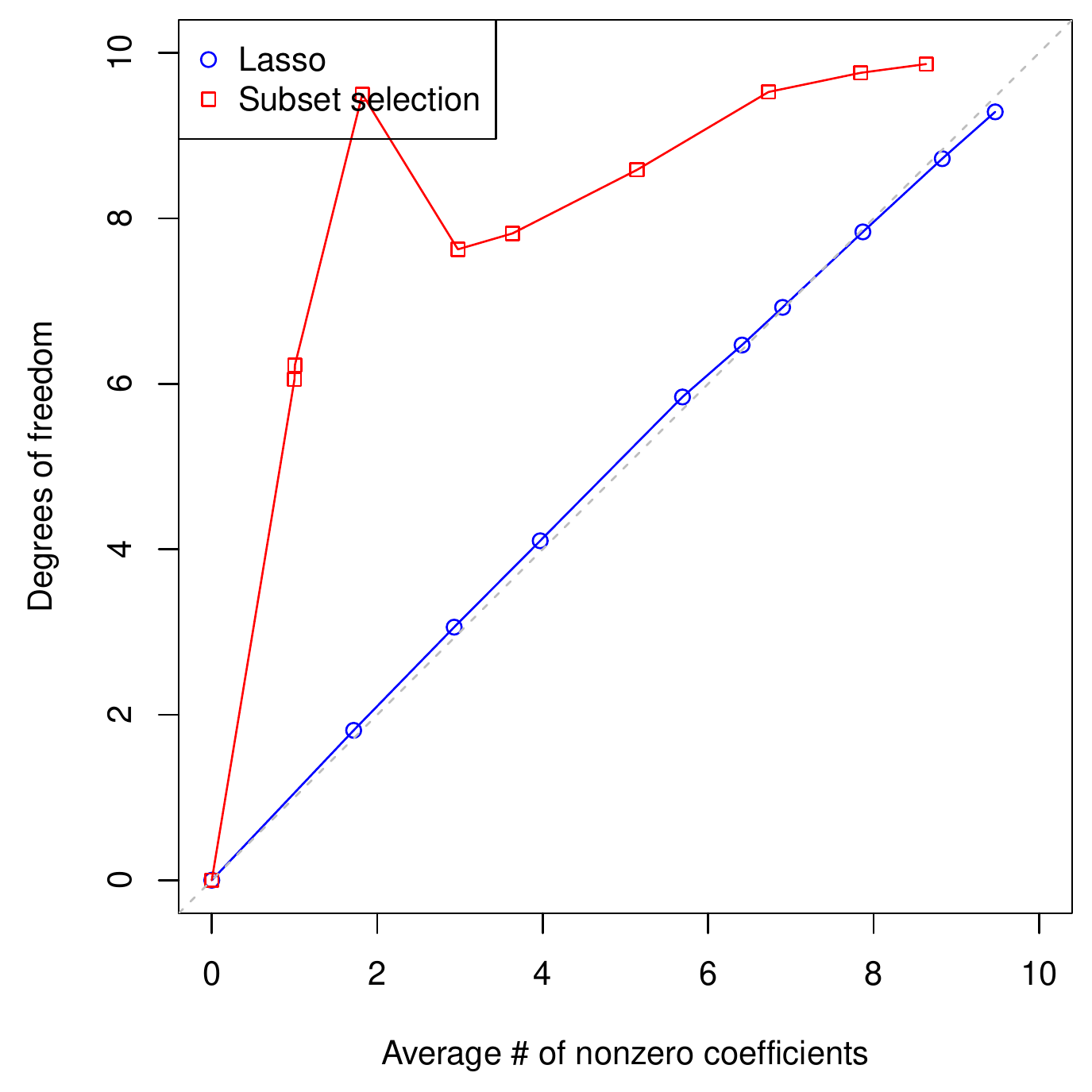}
\caption{\small\it A simulated regression example with $n=20$,
$p=10$.  We drew 100 copies of the outcome $y$ from the same
sparse regression setup, and fit the lasso and best subset selection
estimates each time, across 10 prespecified tuning parameter values.
The plot shows the average number of selected variables by the lasso
(in blue) and best subset selection (in red), across the tuning
parameter values, versus their (estimated) degrees of freedom.  The
lasso degrees of freedom lines up with the number of
selected variables, but the same is not true for subset selection,
with its degrees of freedom being relatively much larger.} 
\label{fig:lassosubset}
\end{figure}

Figure \ref{fig:lassosubset} shows an empirical comparison between
the degrees of freedom of the lasso and best subset selection fits,
for a simple example with $n=20$, $p=10$.  The predictor
variables were setup to have a block correlation structure,
in that variables 1 through 4 had high pairwise correlation (between
0.6 and 0.9), variables 5 through 10 also had high pairwise
correlation (between 0.6 and 0.9), and the two blocks were
uncorrelated with each other.  The outcome $y$ was drawn by adding
independent normal noise to $X\beta^*$, for some true coefficient
vector $\beta^*$, supported on the first block of variables, and on
one variable in the second block. 
We computed the lasso estimate in \eqref{eq:lasso} over 10 values of
the tuning parameter $\lambda$, as well as a best subset selection 
estimate  
\begin{equation}
\label{eq:subset}
\hbetasub \in \argmin_{\beta \in \R^p} \,  
\half\|y-X\beta\|_2^2 + \lambda\|\beta\|_0,
\end{equation}
over its own 10 values of $\lambda$.  Recall that $\|\beta\|_0 =  
\sum_{j=1}^p 1\{\beta_j \not= 0\}$. We repeated this process 100
times, i.e., drew 100 copies of $y$ from the
described regression model, keeping $X$ and $\beta^*$ fixed, and each 
time computed fitted values from the lasso and best subset
selection across the same 10 values of the tuning parameter.  For
each method and value of $\lambda$, we then:
\begin{enumerate}
\item computed the average number of nonzero coefficients over the 100 
  trials;
\item evaluated the covariance in \eqref{eq:df} empirically across the
  100 trials, as an (unbiased) estimate of the degrees of freedom. 
\end{enumerate}
Figure \ref{fig:lassosubset} plots the first quantity
versus the second quantity, with the lasso in blue and best subset
selection in red.  As prescribed by Theorem \ref{thm:lassodf}, the
(estimated) degrees of freedom of the lasso fit is closely aligned with
the average number of nonzero coefficients in its estimate.  But
subset selection does not follow the same trend; its (estimated)
degrees of freedom is much larger than its delivered number of nonzero  
coefficients.  For example, when $\lambda$ is tuned so that the 
subset selection estimate has a little less than 3 nonzero 
coefficients on average, the fit uses about 9 degrees of freedom.   

Why does this happen?  Again, this can be intuitively explained by
shrinkage---this time, a lack thereof.  If we denote the support
of a best subset selection solution by
\smash{$\cAsub=\mathrm{supp}(\hbetasub)$}, and abbreviate
$\cA=\cAsub$, then it is not hard to see that 
\begin{equation*}
\hbetasub_\cA = (X_\cA^T X_\cA)^{-1} X_\cA^T y,
\end{equation*}
i.e., the active coefficients are given by least squares on
the active variables $X_\cA$ (the submatrix of $X$ formed by 
taking the columns in $\cA$).  Therefore, like the lasso, best subset 
selection chooses an active set of variables adaptively, but unlike
the lasso, it fits their coefficients without shrinkage, using
ordinary least squares. It pays for the ``surplus'' of covariance
from the adaptive model search, as well as the usual amount from least
squares estimation, resulting in a total degrees of freedom much larger
than $|\cA|$ (or rather, $\E|\cA|$).

A clarifying note: simulations along the lines of that in Figure
\ref{fig:lassosubset} can be found throughout the literature and we do
not mean to claim originality here (e.g., see Figure 4 of \citet{covinf}
for an early example, and Figure 2 of \citet{dfflawed} for a recent
example).  This simulation is instead simply meant to motivate the
work that follows, as an aim of this paper is to examine the observed
phenomenon in Figure \ref{fig:lassosubset} more formally.

\subsection{Degrees of freedom and optimism}
\label{sec:opt}

Degrees of freedom is closely connected
to the concept of optimism, and so alternatively, we could have
motivated the study of the covariance term on the right-hand side in
\eqref{eq:df} from the perspective of the optimism, rather than the
complexity, of a fitting procedure. 
Assuming only that $y$ is drawn from the model in
\eqref{eq:model}, and that $y'$ is an independent copy of $y$ (i.e., an
independent draw from \eqref{eq:model}), it is straightforward to show
that for any fitting procedure $f$,
\begin{equation}
\label{eq:opt}
\E\|y' - f(y)\|_2^2  - \E\|y - f(y)\|_2^2  = 2 \sigma^2 \cdot \df(f).
\end{equation}
The quantity on the left-hand side above is called the {\it optimism}
of $f$, i.e., the difference in the mean squared test error and mean
squared training error.  The identity in \eqref{eq:opt} shows that
(for uncorrelated, homoskedastic regression errors as in
\eqref{eq:model}) the optimism of $f$ is just a positive constant
times its degrees of freedom; in other words, fitting procedures
with a higher degrees of freedom will have higher a optimism.  Hence,  
from the example in the last section, we know when they are tuned to
have the same (expected) number of variables in the
fitted model, best subset selection will produce a training error that is 
generally far more optimistic than that produced by the lasso. 

\subsection{Lagrange versus constrained problem forms}

Recall that we defined the subset selection estimator using the
Lagrange form optimization problem \eqref{eq:subset}, instead of 
the (perhaps more typical) constrained form definition  
\begin{equation}
\label{eq:constr}
\hbetasub \in \argmin_{\beta \in \R^p} \,  
\|y-X\beta\|_2^2 \;\,\st\;\, \|\beta\|_0 \leq k.
\end{equation}
There are several points now worth making.  First, these are nonconvex
optimization problems, and so the Lagrange and constrained forms
\eqref{eq:subset} and \eqref{eq:constr} of subset selection are
generally not equivalent.  In fact, for all $\lambda$, solutions of
\eqref{eq:subset} are solutions of \eqref{eq:constr} for some choice
of $k$, but the reverse is not true.  Second, even in
situations in which the Lagrange and constrained forms of a particular
optimization problem are equivalent (e.g., this is true under strong
duality, and 
so it is true for most convex problems, under very weak conditions),
there is a difference between studying the degrees of freedom of
an estimator defined in one problem form versus the other.  This is 
because the map from the Lagrange parameter in one form to the
constraint bound in the other generically depends on $y$, i.e., it
is a random mapping (\citet{morereglessdf} discuss this for ridge
regression and the lasso).

Lastly, in this paper, we focus on the Lagrange form \eqref{eq:subset} of
subset selection because we find this problem is easier to analyze 
mathematically.  For example, in Lagrange form with $X=I$,  
the $i$th component of the subset selection fit \smash{$\hbetasub_i$}
depends on $y_i$ only (and is given by hard thresholding), for each 
$i=1,\ldots n$; in constrained form with $X=I$, each
\smash{$\hbetasub_i$} is a function of the order statistics of
$|y_1|,\ldots |y_n|$, and hence depends on the whole sample.  

Given the general spirit of our paper, it is important to recall
the relevant work of \citet{ye}, who studied degrees of freedom for
special cases of best subset selection in constrained form.  In one
such special case (orthogonal predictors with null underlying 
signal), the author derived a simple expression for degrees of
freedom as the sum of the $k$ largest order statistics from a sample
of $n$ independent $\chi_1^2$ random variables.  This indeed
establishes that, in this particular special case, the constrained
form of best subset selection with $k$ active variables has degrees of
freedom larger than $k$.   It does not, however, imply any results
about the Lagrange case for the reasons explained above.

\subsection{Assumptions, notation, and outline}

Throughout this work, we will assume the model
\begin{equation}
\label{eq:normal}
y = \mu + \epsilon, \;\;\; \epsilon \sim N(0,\sigma^2 I). 
\end{equation}
Note that this is stronger than the model in \eqref{eq:model}, since
we are assuming a normal error distribution.  While the
model in \eqref{eq:model} is sufficient to define the notion of
degrees of freedom in general, 
we actually require normality for the calculations to
come---specifically, Lemma \ref{lem:htdf} (on the degrees of freedom
of hard thresholding), and all results in Section \ref{sec:steinex} 
(on extending Stein's formula), rely on the normal error model.
Beyond this running assumption, we will make any additional
assumptions clear when needed.

In terms of notation, we write $M^+$ to denote the (Moore-Penrose)
pseudoinverse of a matrix $M$, with $M^+=(M^T M)^+ M^T$ for
rectangular matrices $M$, and we write $M_S$ to denote the submatrix
of $M$ whose columns correspond to the set of indices $S$.  We write
$\phi$ for the standard normal density function and $\Phi$ for the
standard normal cumulative distribution function.

Finally, here is an outline for the rest of this article.  In
Section \ref{sec:orthx}, we derive an explicit formula for the
degrees of freedom of the best subset selection fit, under
orthogonal predictors $X$.  We also introduce the notion of
search degrees of freedom for subset selection, and study its
characteristics in various settings.  In Section \ref{sec:sdf}, we
define search degrees of freedom for generic adaptive regression
procedures, including the lasso and ridge regression as special
cases. Section \ref{sec:genx} returns to considering best subset 
selection, this time with general predictor variables $X$. Because 
exact formulae for the degrees of freedom and search degrees of
freedom of best subset selection are not available in the general $X$ 
case, we turn to simulation to investigate these quantities. 
We also examine the search degrees of freedom of the lasso across the
same simulated examples (as its analytic calculation is again
intractable for general $X$).  Section \ref{sec:steinex}
casts all of this work on degrees of freedom (and search degrees of 
freedom) in a different light, by deriving an extension of Stein's
formula.  Stein's formula is a powerful tool that can be used to
compute the degrees of freedom of continuous and almost differentiable
fitting procedures; our extension covers functions that have
``well-behaved'' points of discontinuity, in some sense.  This
extended version of Stein's formula offers an alternative proof of the
exact result in Section \ref{sec:orthx} (the orthogonal $X$ case), and  
potentially, provides a perspective from which we can formally
understand the empirical findings in Section \ref{sec:genx} (the
general $X$ case).  In Section \ref{sec:discuss}, we conclude with
some discussion.

\section{Best subset selection with an orthogonal $X$} 
\label{sec:orthx}

In the special case that $X \in \R^{n\times p}$ is orthogonal, i.e., $X$ has 
orthonormal columns, we can compute the degrees of freedom of the best
subset selection fit directly.

\begin{theorem}
\label{thm:subsetdf}
Assume that $y \sim N(\mu,\sigma^2 I)$, and
that $X$ is orthogonal, meaning that $X^T X = I$.  Then the best  
subset selection fit \smash{$\hmusub=X\hbetasub$}, at any
fixed value of $\lambda$, has degrees of freedom  
\begin{equation}
\label{eq:subsetdf}
\df(\hmusub) = \E|\cAsub| +  
\frac{\sqrt{2\lambda}}{\sigma}
\sum_{i=1}^p
\bigg[\phi\bigg(\frac{\sqrt{2\lambda}-(X\T\mu)_i}{\sigma}\bigg) +  
\phi\bigg(\frac{\sqrt{2\lambda}+(X\T\mu)_i}{\sigma}\bigg)\bigg],
\end{equation}
where $\phi$ is the standard normal density.
\end{theorem}

The proof essentially reduces to a calculation on the degrees of
freedom of the (componentwise) hard thresholding operator because, in 
the orthogonal $X$ case, the best subset selection solution is exactly 
hard thresholding of $X^T y$.  Formally, define the hard
thresholding operator $H_t : \R^n \rightarrow \R^n$, at a fixed level
$t \geq 0$, by its coordinate functions 
\begin{equation*}
[H_t(y)]_i = y_i \cdot 1\{|y_i| \geq t\}, \;\;\; i=1,\ldots n.
\end{equation*}
Let $\cA_t$ denote the support set of the output,
$\cA_t=\supp(H_t(y))$.  The following result simply comes from the
normality of $y$, and the definition of degrees of freedom in 
\eqref{eq:df}.  

\begin{lemma}
\label{lem:htdf}
Assume that $y \sim N(\mu,\sigma^2 I)$, and $t \geq 0$ is arbitrary
but fixed.  Then the hard thresholding operator $H_t$ has degrees of 
freedom  
\begin{equation}
\label{eq:htdf}
\df(H_t) = \E|\cA_t| + \frac{t}{\sigma}
\sum_{i=1}^n \bigg[\phi\bigg(\frac{t-\mu_i}{\sigma}\bigg) +  
\phi\bigg(\frac{t+\mu_i}{\sigma}\bigg)\bigg].
\end{equation}
\end{lemma}

\noindent
{\it Remark.} This result, on the degrees of freedom of the hard
thresholding operator, can be found in both \citet{sparsenet} and 
\citet{score}.  The former work uses degrees of freedom as a
calibration tool in nonconvex sparse regression; the latter
derives an estimate of the right hand side in \eqref{eq:htdf}
that, although biased, is consistent under some conditions.

\smallskip
\smallskip

The proofs of Lemma \ref{lem:htdf}, and subsequently Theorem
\ref{thm:subsetdf}, involve straightforward calculations, but are
deferred until the appendix for the sake of readability.  

We have established that, for and orthogonal $X$, the
degrees of freedom of the subset selection fit is equal to
$\E|\cAsub|$, plus an ``extra'' term.  We make some observations 
about this term in the next section.    

\subsection{Search degrees of freedom}

The quantity 
\begin{equation}
\label{eq:subsetsdf}
\sdf(\hmusub) = 
\frac{\sqrt{2\lambda}}{\sigma}
\sum_{i=1}^p
\bigg[\phi\bigg(\frac{\sqrt{2\lambda}-(X\T\mu)_i}{\sigma}\bigg) +  
\phi\bigg(\frac{\sqrt{2\lambda}+(X\T\mu)_i}{\sigma}\bigg)\bigg]
\end{equation}
appearing in \eqref{eq:subsetdf} is the amount by which the degrees of
freedom exceeds the expected number of selected variables.
We will refer to this the {\it search degrees of freedom} of best
subset selection, because roughly speaking, we can think of it as the
extra amount of covariance that comes from searching through the space 
of models. Note that \smash{$\sdf(\hmusub) > 0$} for any 
$\lambda>0$, because the normal density is supported everywhere, and
therefore we can indeed conclude that 
\smash{$\df(\hmusub) > \E|\cAsub|$}, as we suspected, in the case of
an orthogonal predictor matrix.   

How big is $\sdf(\hmusub)$? At the extremes:
$\sdf(\hmusub)=0$ when $\lambda=0$, and
$\sdf(\hmusub) \rightarrow 0$ as
$\lambda\rightarrow\infty$. In words, searching has no cost
when all of the variables, or none of the variables, are in 
the model.  But the behavior is more interesting for intermediate
values of $\lambda$.  The precise shape of the search degrees of
freedom curve \eqref{eq:subsetsdf}, over $\lambda$, depends on the
underlying signal $\mu$; the next three sections study three canonical
cases for the underlying signal. 

\subsection{Example: null signal}

\begin{figure}[tb]
\centering
\includegraphics[width=0.45\textwidth]{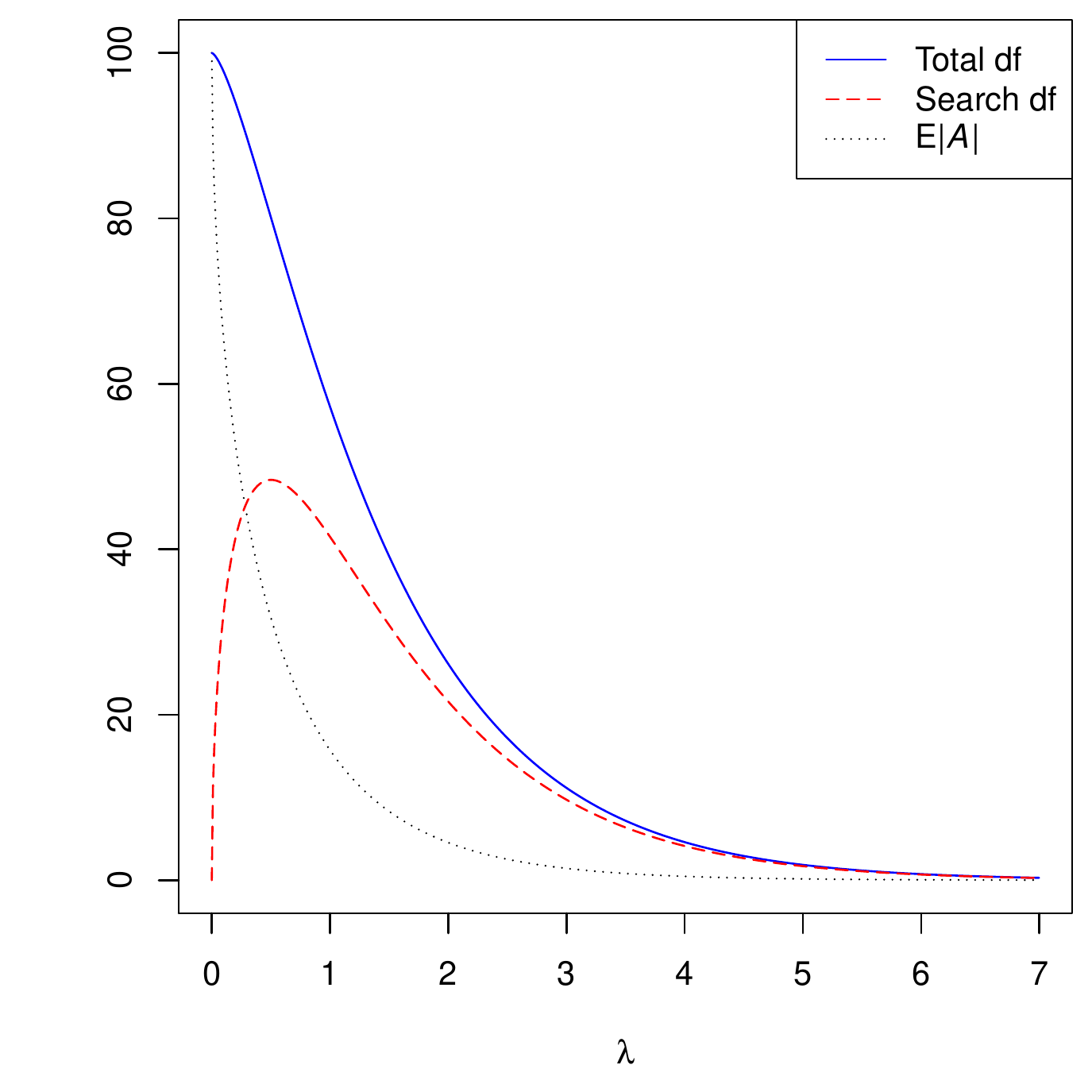}
\hspace{5pt}
\includegraphics[width=0.45\textwidth]{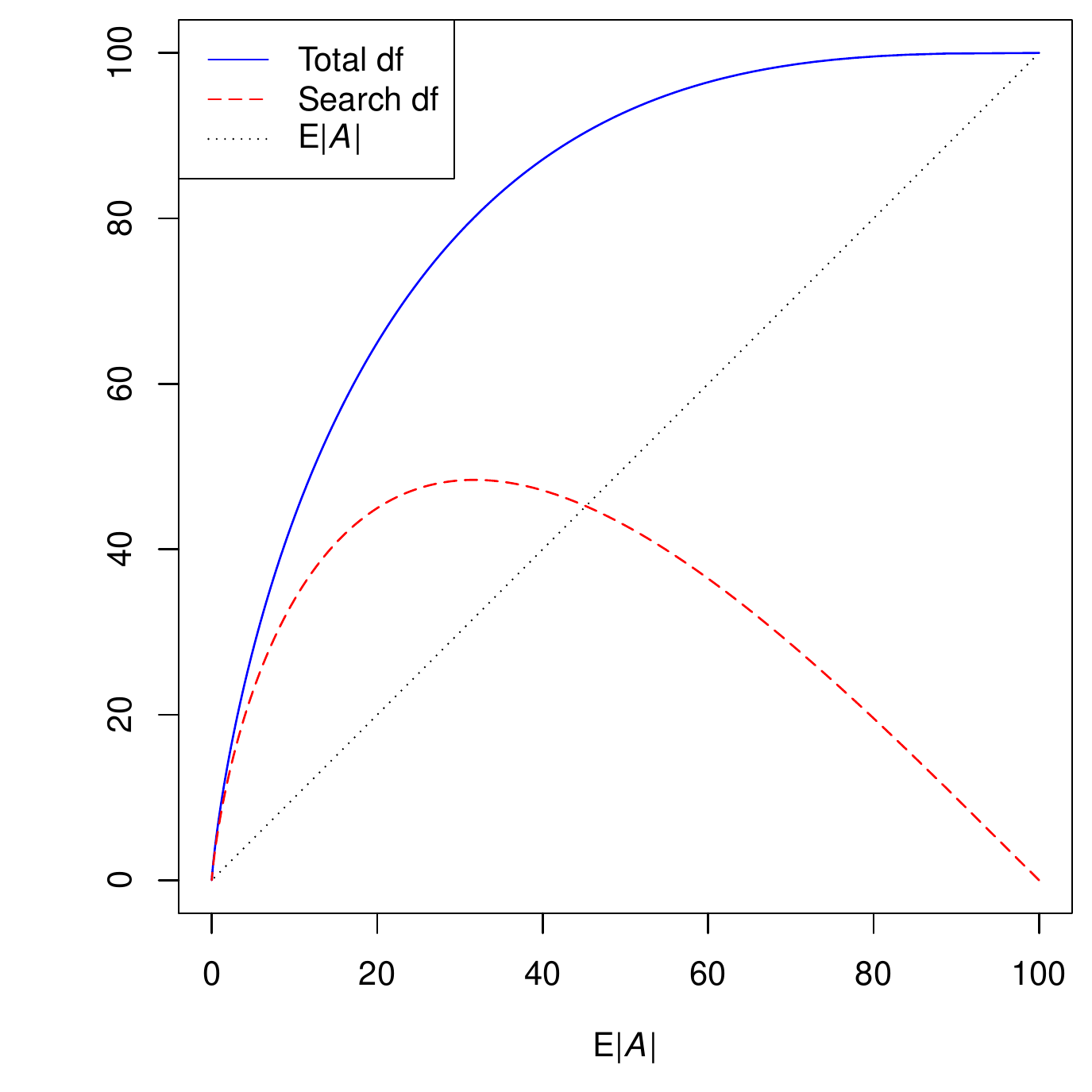}
\caption{\small\it An example with $n=p=100$, $X=I$, and $\mu=0$. 
  The left panel plots the curves \smash{$\df(\hmusub)$},
  \smash{$\sdf(\hmusub)$}, and $\E|\cAsub|$ as functions of $\lambda$,
  drawn as blue, red, and black lines, respectively. The right panel
  plots the same quantities with respect to $\E|\cAsub|$.}
\label{fig:subsetdfnull}
\end{figure}

We consider first the case of a null underlying signal, i.e., $\mu=0$. 
The best subset selection search degrees of freedom
\eqref{eq:subsetsdf}, as a function of $\lambda$, becomes  
\begin{equation}
\label{eq:subsetsdfnull}
\sdf(\hmusub) = 
\frac{2p\sqrt{2\lambda}}{\sigma} 
\phi\bigg(\frac{\sqrt{2\lambda}}{\sigma}\bigg).
\end{equation}
In Figure \ref{fig:subsetdfnull}, we plot the quantities
$\df(\hmusub)$, $\sdf(\hmusub)$, and $\E|\cAsub|$ as functions of
$\lambda$, for a simple example with $n=p=100$, underlying signal
$\mu=0$, noise variance $\sigma^2=1$, and predictor matrix 
$X=I$, the $100 \times 100$ identity matrix.  We emphasize that this
figure was produced without any random draws or simulations, and
the plotted curves are exactly as prescribed by Theorem
\ref{thm:subsetdf}  (recall that $\E|\cAsub|$ also has an explicit
form in terms of $\lambda$, given in the proof of Lemma
\ref{lem:htdf}).  In the left panel, we can see that the search
degrees of freedom curve is maximized at approximately $\lambda =
0.5$, and achieves a maximum value of nearly 50. 
That is, when $\lambda = 0.5$, best subset selection spends
nearly 50 (extra) parameters searching through the space of models!   

It is perhaps more natural to parametrize the curves in terms
of the expected number of active variables $\E|\cAsub|$ (instead of 
$\lambda$), as displayed in the right panel of Figure
\ref{fig:subsetdfnull}. This parametrization reveals something 
interesting: the search degrees of freedom curve is maximized at
roughly $\E|\cAsub|=31.7$. In other words, 
searching is most costly when there are approximately 31.7 variables
in the model. This is a bit counterintuitive, because there are
more subsets of size $50$ than any other size, that is, the function  
\begin{equation*} 
F(k) = {p \choose k}, \; k=1,2,\ldots p,
\end{equation*} 
is maximized at $k=p/2=50$. Hence we might believe that
searching through subsets of variables is most costly when
$\E|\cAsub|=50$, because in this case the search space is largest. 
Instead, the maximum actually occurs at about $\E|\cAsub|=31.7$.
Given the simple form \eqref{eq:subsetsdfnull} of the search degrees 
of freedom curve in the null signal case, we can verify this 
observation analytically: direct calculation shows that the right hand
side in \eqref{eq:subsetsdfnull} is maximized at $\lambda =
\sigma^2/2$, which, when plugged into the formula for the expected
number of selected variables in the null case,
\begin{equation*}
\E|\cAsub| = 2p \Phi\bigg(\frac{-\sqrt{2\lambda}}{\sigma}\bigg),
\end{equation*}
yields $\E|\cAsub| = 2 \Phi(-1) p \approx 0.317 p$.

Although this calculation may have been reassuring, the intuitive
question remains: why is the 31.7 variable model associated with
the highest cost of model searching (over, say, the 50 variable
model)?  At this point, we cannot offer a truly satisfying
intuitive answer, but we will attempt an explanation nonetheless.
Recall that search degrees
of freedom measures the additional amount of covariance in
\eqref{eq:df} that we attribute to searching through the space of
models---additional from the baseline amount $\E|\cAsub|$, which  
comes from estimating the coefficients in the selected model.
The shape of the search degrees of freedom curve, when $\mu=0$,
tells us that there is more covariance to be gained when the selected
model has 31.7 variables than when it has 50 variables.  
As the size of the selected subset $k$ increases from 0 to 50, note
that: 
\begin{enumerate}
\item the number of subsets of size $k$
  increases, which means that there are more opportunities to decrease  
  the training error, and so the total degrees of freedom (optimism) 
  increases; 
\item trivially, the baseline amount of fitting also
  increases, as this baseline is just $k$, the degrees of freedom  
  (optimism) of a fixed model on $k$ variables. 
\end{enumerate}
Search degrees of freedom is the difference between these two
quantities (i.e., total minus baseline degrees of freedom), and as it
turns out, the two are optimally balanced at approximately $k=31.7$
(at exactly $k=2\Phi(-1)p$) in the null signal case.

\subsection{Example: sparse signal}

\begin{figure}[tb]
\centering
\includegraphics[width=0.45\textwidth]{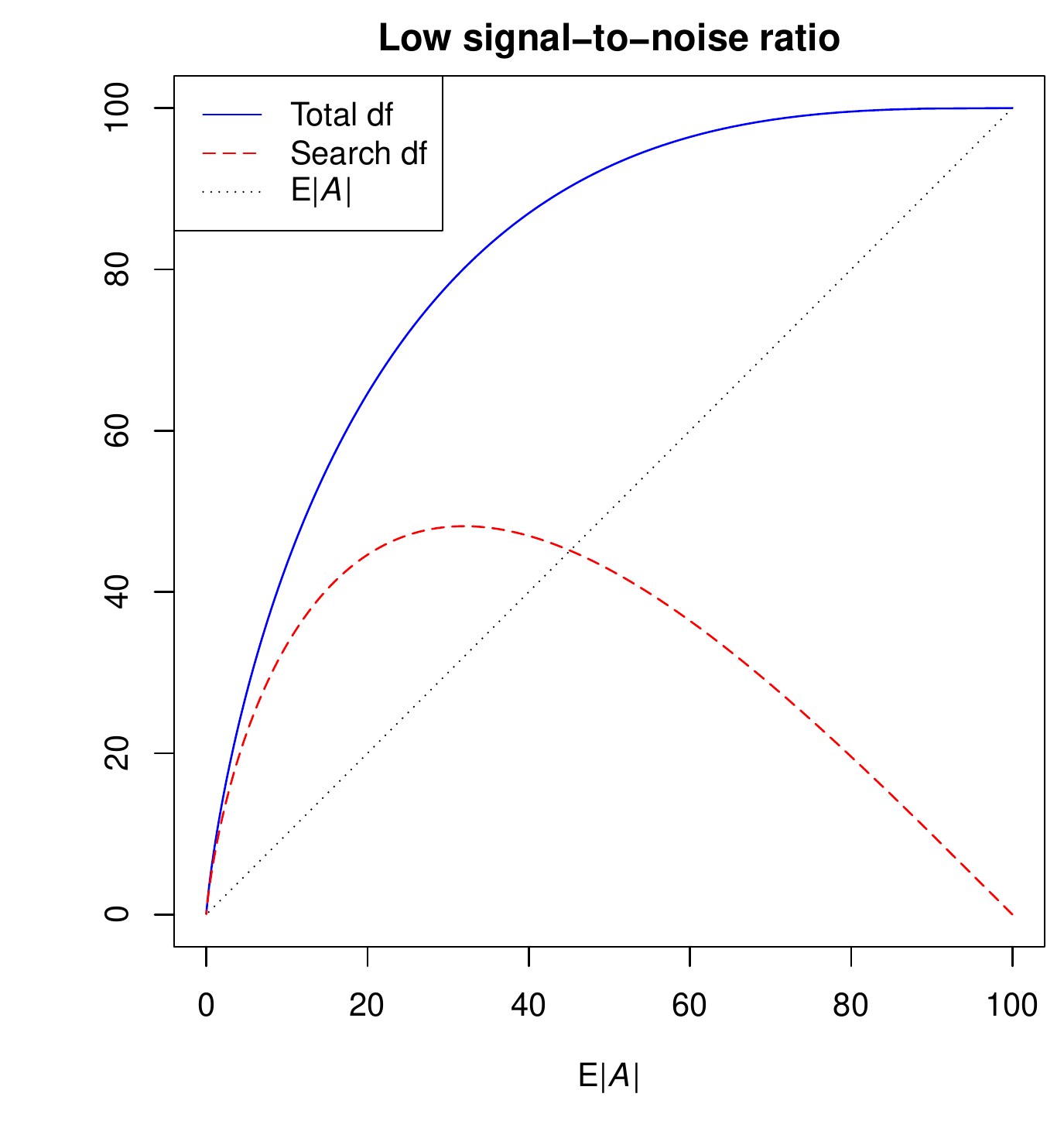}
\hspace{5pt}
\includegraphics[width=0.45\textwidth]{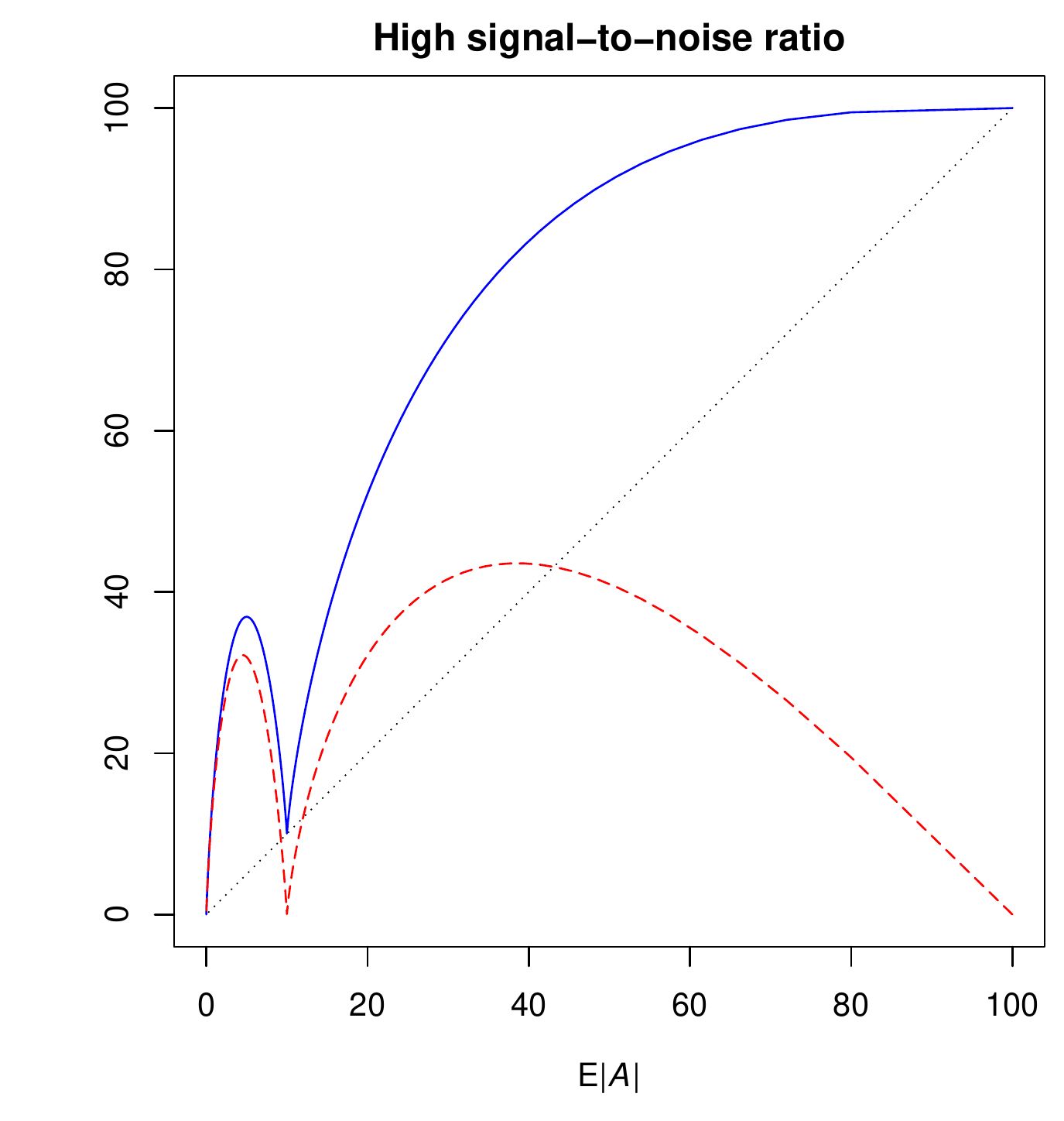}
\caption{\small\it  An example with $n=p=100$, $X=I$, and
  $\mu=X\beta^*$ with $\beta^*$ as in \eqref{eq:sparsecoef}.  The left
  panel corresponds to the choice $\rho=1$ (low signal-to-noise
  regime) and the right to $\rho=8$ (high signal-to-noise regime).} 
\label{fig:subsetdfsparse}
\end{figure}

Now we consider the case in which $\mu=X\beta^*$, for some sparse 
coefficient vector $\beta^* \in \R^p$.  We let $\cA^* = \supp(\beta^*)$
denote the true support set, and $k^*=|\cA^*|$ the true number of
nonzero coefficients, assumed to be small.  The search degrees of
freedom curve in \eqref{eq:subsetsdf} is 
\begin{equation}
\label{eq:subsetsdfsparse}
\sdf(\hmusub) = 
\frac{\sqrt{2\lambda}}{\sigma}
\sum_{i \in \cA^*}
\bigg[\phi\bigg(\frac{\sqrt{2\lambda}-\beta^*_i}{\sigma}\bigg) +  
\phi\bigg(\frac{\sqrt{2\lambda}+\beta^*_i}{\sigma}\bigg)\bigg] 
+ \frac{2(p-k^*)\sqrt{2\lambda}}{\sigma}
\phi\bigg(\frac{\sqrt{2\lambda}}{\sigma}\bigg).
\end{equation}
When the nonzero coefficients $\beta^*_i$ are moderate (not
very large), the curve in \eqref{eq:subsetsdfsparse} acts much like the
search degrees of freedom curve \eqref{eq:subsetsdfnull} in the null
case.  Otherwise, it can behave very differently.  We therefore
examine two different sparse setups by example, having low and high 
signal-to-noise ratios.  See Figure \ref{fig:subsetdfsparse}.  In
both setups, we take $n=p=100$, $\sigma^2=1$, $X=I$, and
$\mu=X\beta^*$, with
\begin{equation}
\label{eq:sparsecoef}
\beta^*_i = \begin{cases}
\rho & i=1,\ldots 10 \\
0 & i=11,\ldots 100.
\end{cases} 
\end{equation}
The left panel uses $\rho=1$, and the right uses $\rho=8$.  We plot
the total degrees of freedom and search degrees of freedom of 
subset selection as a function of the expected number of selected
variables (and note, as before, that these plots are produced by
mathematical formulae, not by simulation). The curves in the left
panel, i.e., in the low signal-to-noise ratio case, appear extremely
similar to those in 
the null signal case (right panel of Figure \ref{fig:subsetdfnull}).
The search degrees of freedom curve peaks when the expected
number of selected variables is about $\E|\cAsub|=31.9$, and its
peak height is again just short of 50.  

Meanwhile, in the high
signal-to-noise ratio case, i.e., the right panel of Figure
\ref{fig:subsetdfsparse}, the behavior is very different.  
The search degrees of freedom curve is bimodal, and is basically  
zero when the expected number of selected variables is 10.
The intuition: with such a high signal-to-noise ratio in the true
model \eqref{eq:sparsecoef}, best subset selection is able to select
the same (true) subset of 10 variables in every random data instance,
and therefore the size 10 model produced by subset selection is akin
to a fixed model, with no real searching performed whatsoever.
Another interesting point is that the cost of model searching is
very high when the selected model has average size equal to 5; here
the search component contributes over 30 degrees of freedom to the
total.  Intuitively, with 10 strong variables in the true model
\eqref{eq:sparsecoef}, there are many competitive subsets of size 5,
and hence a lot of searching is needed in order to report the best
subset of size 5 (in terms of training error). 

It is worth mentioning the interesting, recent works of
\citet{morereglessdf} and \citet{dfflawed}, which investigate
unexpected nonmonoticities in the (total) degrees of freedom of an
estimator, as a function of some underlying parametrization for the
amount of imposed regularization. We note that the right panel of
Figure \ref{fig:subsetdfsparse} portrays a definitive example of this, 
in that the best subset selection degrees of freedom undergoes a
major nonmonoticity at 10 (expected) active variables, as discussed  
above. 

\subsection{Example: dense signal}

\begin{figure}[tb]
\centering
\includegraphics[width=0.45\textwidth]{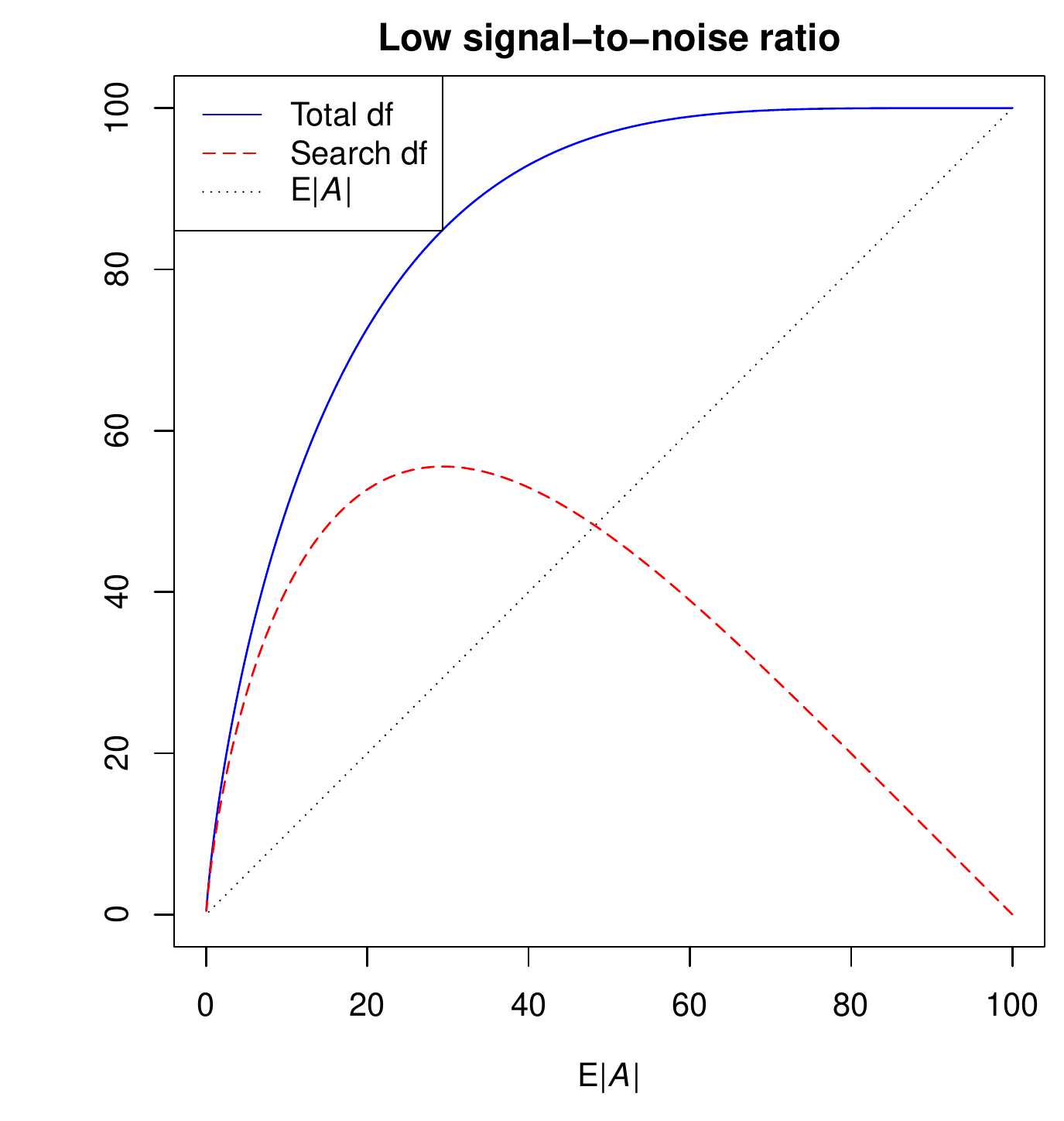}
\hspace{5pt}
\includegraphics[width=0.45\textwidth]{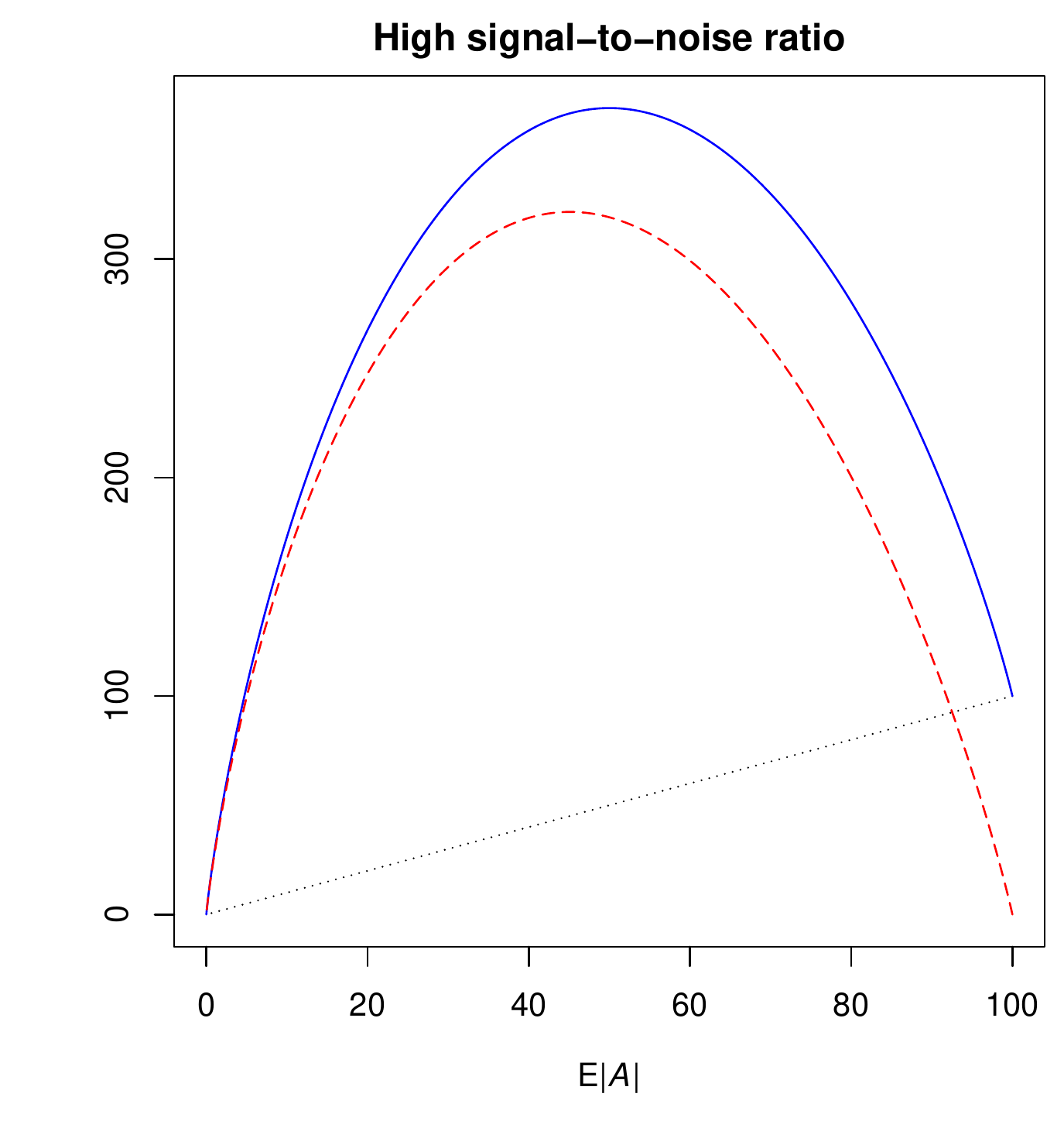}
\caption{\small\it  An example with $n=p=100$, $X=I$, and
  $\mu=X\beta^*$ with $\beta^*_i=\rho$, $i=1,\ldots p$.
  The left panel corresponds to $\rho=1$ (low signal-to-noise regime)
  and the right to $\rho=8$ (high signal-to-noise regime).} 
\label{fig:subsetdfdense}
\end{figure}

The last case we consider is that of a dense underlying signal, 
$\mu=X\beta^*$ for some dense coefficient vector
$\beta^* \in \R^p$.  For the sake of completeness, in the present
case, the expression \eqref{eq:subsetsdf} for the search degrees of
freedom of best subset selection is
\begin{equation}
\label{eq:subsetsdfdense}
\sdf(\hmusub) = 
\frac{\sqrt{2\lambda}}{\sigma}
\sum_{i=1}^p
\bigg[\phi\bigg(\frac{\sqrt{2\lambda}-\beta^*_i}{\sigma}\bigg) +  
\phi\bigg(\frac{\sqrt{2\lambda}+\beta^*_i}{\sigma}\bigg)\bigg].
\end{equation}
The search curve degrees of freedom curve \eqref{eq:subsetsdfdense}
exhibits a very similar behavior to the curve \eqref{eq:subsetsdfnull} in
the null signal case when 
the coefficients $\beta^*_i$ are small or moderate, but a very
different behavior when some coefficients $\beta^*_i$ are large. 
In Figure \ref{fig:subsetdfdense}, we take $n=p=100$, $X=I$, and $\mu
= X\beta^*$ with $\beta^*_i=\rho$, $i=1,\ldots p$.  The left
panel of the figure corresponds to $\rho=1$, and the right corresponds
to $\rho=8$.  Both panels plot degrees of freedom against the expected
number of selected variables (and, as in the last two subsections, these
degrees of freedom curves are plotted according to their closed-form
expressions, they are not derived from simulation).  We can see that
the low signal-to-noise ratio case, in the left panel,
yields a set of curves quite similar to those from the null signal
case, in the right panel of Figure \ref{fig:subsetdfnull}.  One
difference is that the search degrees of freedom curve has a higher
maximum (its value about 56, versus 48 in the null signal
case), and the location of this maximum is further to the left (occuring
at about $\E|\cAsub|=29.4$, versus $\E|\cAsub|=31.7$ in the former
case).  

On the other hand, the right panel of the figure shows 
the high signal-to-noise ratio case, where the total degrees of
freedom curve is now nonmonotone, and reaches its maximum at an
expected number of selected variables (very nearly) $\E|\cAsub|=50$.
The search degrees of freedom curve itself peaks much later than it
does in the other 
cases, at approximately $\E|\cAsub|=45.2$.  Another striking
difference is the sheer magnitude of the degrees of freedom
curves: at 50 selected variables on average, the total degrees of
freedom of the best subset selection fit is well over 300.
Mathematically, this makes sense, as the search degrees of freedom
curve in \eqref{eq:subsetsdfdense} is increasing in $|\beta^*_i|$.
Furthermore, we can liken the degrees of freedom curves in the right 
panel of Figure \ref{fig:subsetdfdense} to those in a small portion
of the plot in the right panel of Figure \ref{fig:subsetdfsparse},
namely, the portion corresponding to $\E|\cAsub| \leq 10$.  The two  
sets of curves here appear similar in shape.  This is intuitively 
explained by the fact that, in the high signal-to-noise ratio regime,
subset selection over a dense true model is similar to
subset selection over a sparse true model, provided that we constrain
our attention in the latter case to subsets of size less than or
equal to the true model size (since under this constraint, the truly  
irrelevant variables in the sparse model do not play much of a role).  

\section{Search degrees of freedom for general procedures}  
\label{sec:sdf}

Here we extend the notion of search degrees of freedom to general
adaptive regression procedures.  Given an outcome
$y\in\R^n$ and predictors $X \in \R^{n\times p}$, we consider a
fitting procedure $f : \R^n \rightarrow \R^n$ of the form
\begin{equation*}
f(y)=X \hbeta^{(f)},
\end{equation*}
for some estimated coefficients \smash{$\hbeta^{(f)} \in \R^p$}.
Clearly, the lasso and best subset selection
are two examples of such a fitting procedure, with the coefficients 
as in \eqref{eq:lasso} and \eqref{eq:subset}, respectively.
We denote \smash{$\cA^{(f)} = \supp(\hbeta^{(f)})$}, the support set
of the estimated coefficients under $f$.  The overall complexity of
$f$ is measured by its degrees of freedom, as defined in \eqref{eq:df}
(just as it is for all fitting procedures), but we may be also
interested in a degree of complexity associated solely with its model
selection component---i.e., we might ask: how many effective
parameters does $f$ spend in simply selecting the active set
\smash{$\cA^{(f)}$}? 

We propose to address this question by developing a notion of search
degrees of freedom for $f$, in a way that generalizes the notion
considered in the last section specifically for subset selection.
Abbreviating \smash{$\cA=\cA^{(f)}$}, we first
define a modified procedure \smash{$\tilde{f}$} that returns the
least squares fit on the active set $\cA$,
\begin{equation*}
\tilde{f}(y) = P_\cA y.
\end{equation*}
where $P_\cA = X_\cA (X_\cA^T X_\cA)^+ X_\cA^T$ is the projection
onto the span of active predictors $X_\cA$ (note the use of the 
pseudoinverse, as $X_\cA$ need not
have full column rank, depending on the nature of the procedure $f$). 
We now define the search degrees of freedom of $f$ as
\begin{align}
\nonumber
\sdf(f) &= \df(\tilde{f}) - \E[\rank(X_\cA)] \\
\label{eq:sdf}
&= \frac{1}{\sigma^2} \sum_{i=1}^n \Cov\big((P_\cA y)_i,y_i\big) 
- \E[\rank(X_\cA)].
\end{align}
The intuition behind this definition: by construction,
\smash{$\tilde{f}$} and $f$ are identical in their selection of the
active set $\cA$, and only differ in how they 
estimate the nonzero coefficients once $\cA$ has been chosen, with 
$\tilde{f}$ using least squares, and $f$ using a possibly different
mechanism.  If $\cA$ were fixed, then a least
squares fit on $X_\cA$ would use $\E[\rank(X_\cA)]$ degrees of 
freedom, and so it seems reasonable to assign the leftover part,
$\df(\tilde{f})-\E[\rank(X_\cA)]$, as the degrees of 
freedom spent by \smash{$\tilde{f}$} in selecting $\cA$ in
the first place, i.e., the amount spent by $f$ in selecting $\cA$
in the first place. 

It may help to discuss some specific cases.  

\subsection{Best subset selection}

When $f$ is the best subset selection fit, we have
\smash{$\tilde{f}=f$}, i.e., subset selection already performs
least squares on the set of selected variables $\cA$.  Therefore, 
\begin{equation}
\label{eq:sdfls}
\sdf(f) = \df(f) - \E|\cA|,
\end{equation}
where we have also used the fact that $X_\cA$ must have linearly
independent columns with best subset selection (otherwise, we could
strictly decrease the $\ell_0$ penalty in \eqref{eq:subset} while
keeping the squared error loss unchanged).  This matches our
definition \eqref{eq:subsetsdf} of search degrees of freedom for
subset selection in the orthogonal $X$ case---it is the total degrees
of freedom minus the expected number of selected variables, with the
total being explicitly computable for orthogonal predictors, as
we showed in the last section.   

The same expression \eqref{eq:sdfls} holds for any fitting
procedure $f$ that uses least squares to estimate the coefficients in
its selected model, because then \smash{$\tilde{f}=f$}. (Note that,
in full generality, $\E|\cA|$ should be replaced again by 
$\E[\rank(X_\cA)]$ in case $X_\cA$ need not have full column rank.)  
An example of another such procedure is forward stepwise regression.  

\subsection{Ridge regression}

For ridge regression, the active model is 
$\cA=\{1,\ldots p\}$ for any draw of the outcome $y$, which means
that the modified procedure \smash{$\tilde{f}$} is just the full
regression fit on $X$, and  
\begin{equation*}
\sdf(f) = \E[\rank(X)]-\E[\rank(X)] = 0.
\end{equation*}
This is intuitively the correct notion of search degrees of freedom
for ridge regression, since this procedure does not perform any kind of 
variable selection whatsoever.  The same logic carries over to any
procedure $f$ whose active set $\cA$ is almost surely
constant. 

\subsection{The lasso}

The lasso case is an interesting one.  We know from the literature
(Theorem \ref{thm:lassodf}) that the degrees of freedom of the lasso
fit is $\E|\cA|$ (when the predictors are in general position), but
how much of this total can we attribute to model searching?  The
modified
procedure \smash{$\tilde{f}$} that performs least squares on the lasso
active set $\cA$ has been called the {\it relaxed lasso} (see
\citet{relax}, who uses the same term to
refer to a  broader family of debiased lasso estimates). We will
denote the relaxed 
lasso fitted values by \smash{$\hmu^\mathrm{relax}=P_{\cAlas} y$}. 
When $X$ has orthonormal columns, it is not hard to see that the
relaxed lasso fit is given by hard thresholding, just like best subset
selection, but this time with threshold level $t=\lambda$.  The
following result hence holds by the same arguments as those in Section  
\ref{sec:orthx} for subset selection.

\begin{theorem}
\label{thm:relaxdf}
If $y\sim N(\mu,\sigma^2)$, and $X^T X = I$, then the
relaxed lasso fit \smash{$\hmu^\mathrm{relax}=P_{\cAlas} y$}, at a
fixed value $\lambda \geq 0$, has degrees of freedom
\begin{equation*}
\df(\hmu^\mathrm{relax}) =
\E|\cAlas| +  
\frac{\lambda}{\sigma}
\sum_{i=1}^p
\bigg[\phi\bigg(\frac{\lambda-(X\T\mu)_i}{\sigma}\bigg) +  
\phi\bigg(\frac{\lambda+(X\T\mu)_i}{\sigma}\bigg)\bigg].
\end{equation*}
Therefore the lasso has search degrees of freedom
\begin{equation}
\label{eq:lassosdf}
\sdf(\hmulas) = 
\frac{\lambda}{\sigma}
\sum_{i=1}^p
\bigg[\phi\bigg(\frac{\lambda-(X\T\mu)_i}{\sigma}\bigg) +  
\phi\bigg(\frac{\lambda+(X\T\mu)_i}{\sigma}\bigg)\bigg].
\end{equation}
\end{theorem}

The search degrees of freedom formulae \eqref{eq:lassosdf} and
\eqref{eq:subsetsdf} are different as functions of $\lambda$, the
tuning parameter, but this is not a meaningful difference;
when each is parametrized by their respective expected number of
selected variables $\E|\cA|$, the two curves are exactly the same, and 
therefore all examples and figures in Section
\ref{sec:orthx} demonstrating the behavior of the search degrees of 
freedom of best subset selection also apply to the lasso.  In a
sense, this is not a surprise, because for orthogonal predictors both
the lasso and subset selection fits reduce to a sequence of marginal 
considerations (thresholds, in fact), and so their search mechanisms
can be equated. 

But for correlated predictors, we might believe that the search
components associated with the lasso and best subset selection
procedures are actually quite different.  Even though our definition
of search degrees of freedom in \eqref{eq:sdf} is not connected to
computation in any way, the fact that best subset selection
\eqref{eq:subset} is NP-hard for a general $X$ may
seem to suggest (at a very loose level) that it somehow ``searches
more'' than the convex lasso problem \eqref{eq:lasso}.  
For many problem setups, this guess (whether or not properly 
grounded in intuition) appears to be true in simulations,  
as we show next.

\section{Best subset selection with a general $X$}
\label{sec:genx}

\begin{figure}[htb]
\centering
\includegraphics[width=0.55\textwidth]{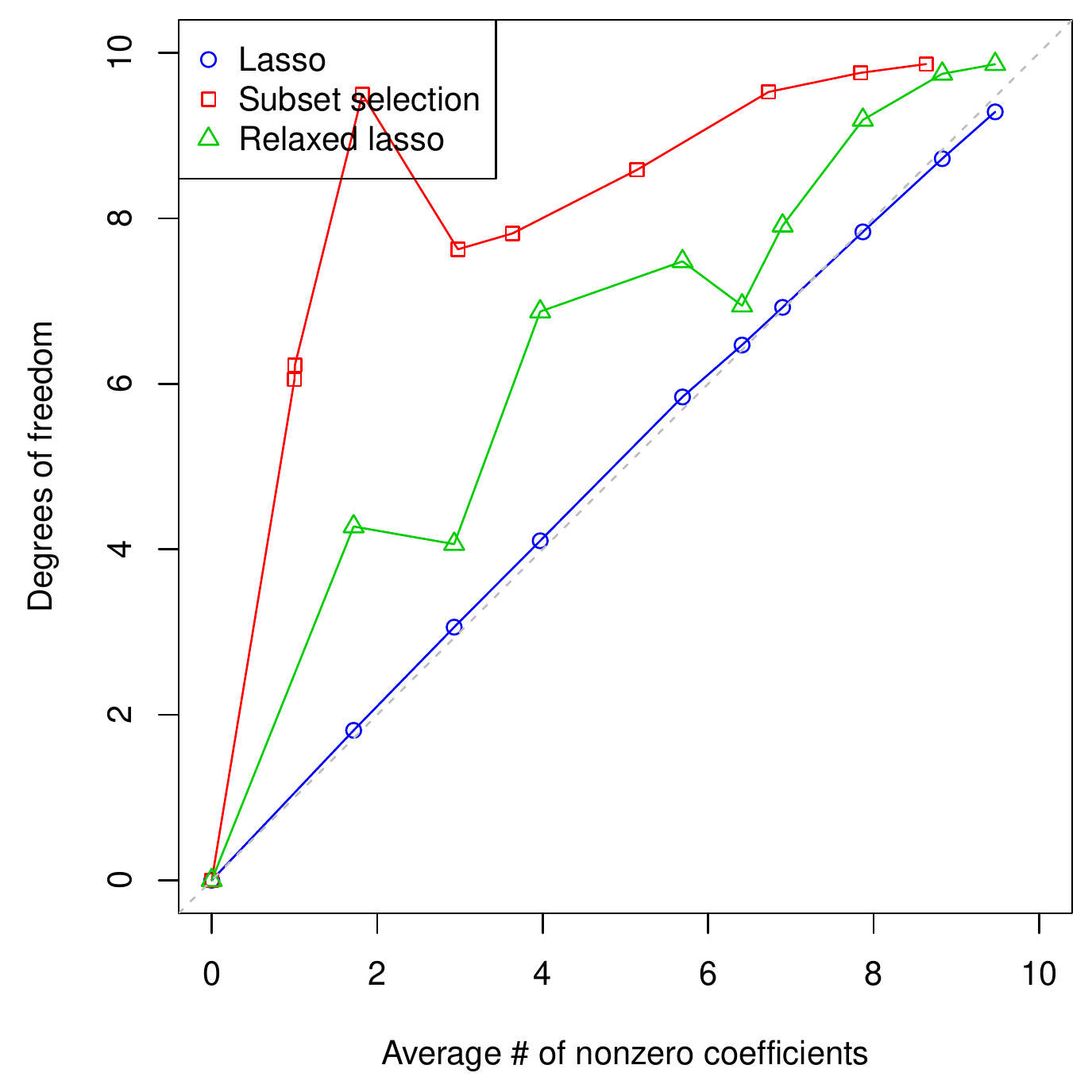}
\caption{\small\it The same simulation setup as in Figure
  \ref{fig:lassosubset}, but now including the relaxed lasso degrees
  of freedom on the left panel, in green.  (The relaxed lasso is the
  fitting procedure that performs least squares on the lasso active
  set.)  We can see that the relaxed lasso has a smaller degrees of
  freedom than subset selection, as a function of their
  (respective) average number of selected variables.  Hence, 
  the lasso exhibits a smaller search degrees of freedom than subset
  selection, in this example.}  
\label{fig:lsr1}
\end{figure}

We look back at the motivating example given in Section
\ref{sec:motivex}, where we estimated the degrees of freedom of 
best subset selection and the lasso by simulation, in a problem
with $n=20$ and $p=10$.  See Section \ref{sec:motivex} for more
details about the setup (i.e., correlation structure of the predictor 
variables $X$, true coefficients $\beta^*$, etc.).  Here we also
consider the degrees of freedom of the relaxed lasso, estimated from
the same set of simulations. Figure \ref{fig:lsr1} plots these
degrees of freedom estimates, in green, on top of
the existing best subset selection and lasso curves from Figure
\ref{fig:lassosubset}. 
Interestingly, the relaxed lasso is seen to have a smaller degrees of
freedom than best subset selection (when each is parametrized by their
own average number of selected variables).  Note that this means the
search degrees of freedom of the lasso (i.e., the difference between
the green curve and the diagonal) is smaller than the search degrees
of freedom of subset selection (the difference between the red
curve and the diagonal). 

This discrepancy between the search degrees of freedom of the lasso
and subset selection, for correlated variables $X$, stands in
contrast to the orthogonal case, where the two quantities were
proven to be equal (subject to the appropriate parametrization). 
Further simulations with correlated predictors show that, for
the most part, this discrepancy persists across a variety of cases;  
consult Figure \ref{fig:lsr2} and the accompanying caption text for
details. However, it is important to note that this phenomenon is not
universal, and in some instances (particularly, when the
computed active set is small, and the true signal is dense) the 
search degrees of freedom of the lasso can grow quite large and
compete with that of subset selection.  Hence, we can see
that the two quantities do not always obey a simple ordering, and the  
simulations presented here call for a more formal understanding
of their relationship.     


\def\d{0.39\textwidth}
\begin{figure}[p]
\centering
\begin{tabular}{m{0.08\textwidth} m{\d} m{\d}}
\hline
Null & 
\includegraphics[width=\d]{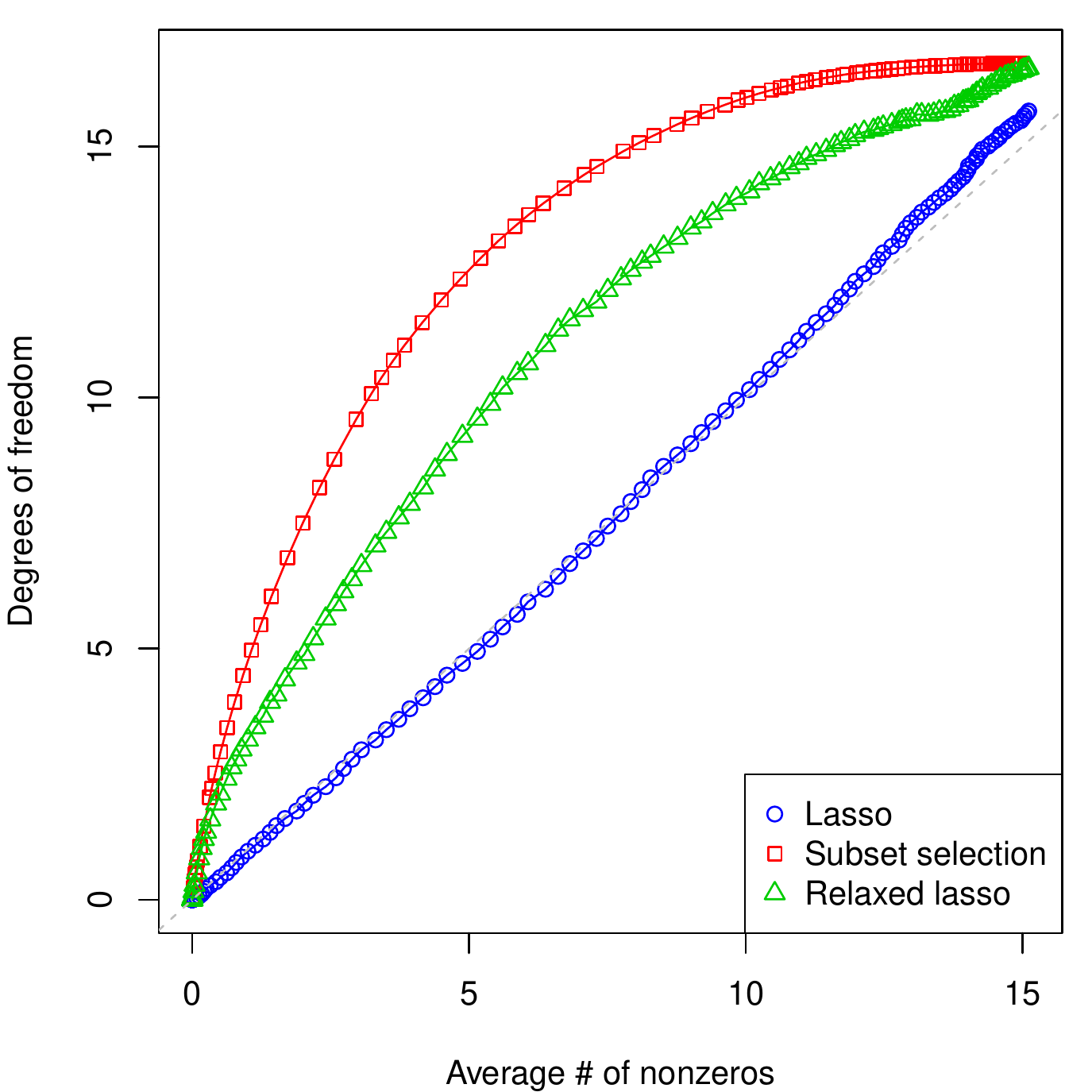} &
\includegraphics[width=\d]{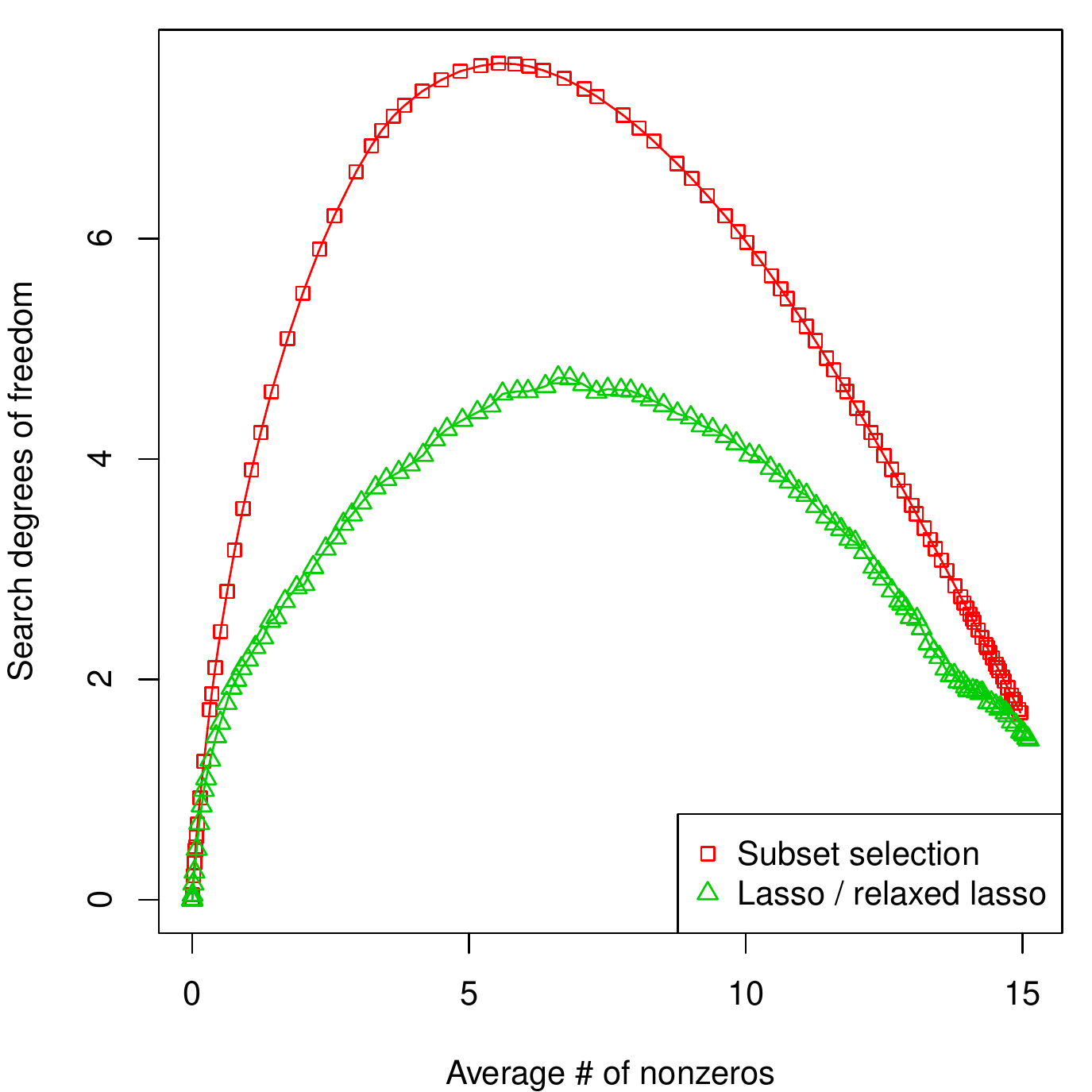} \\
\hline
Sparse &
\includegraphics[width=\d]{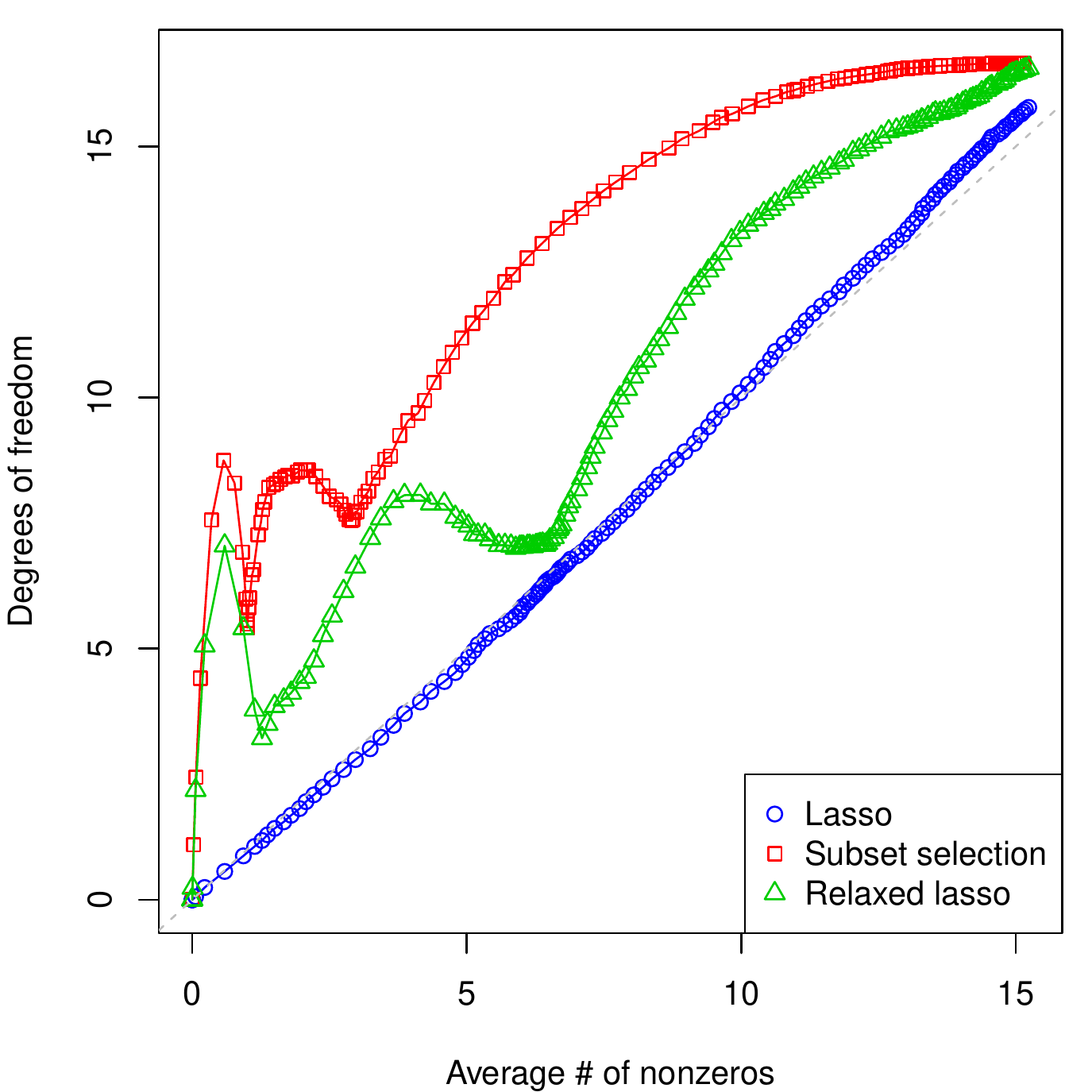} &
\includegraphics[width=\d]{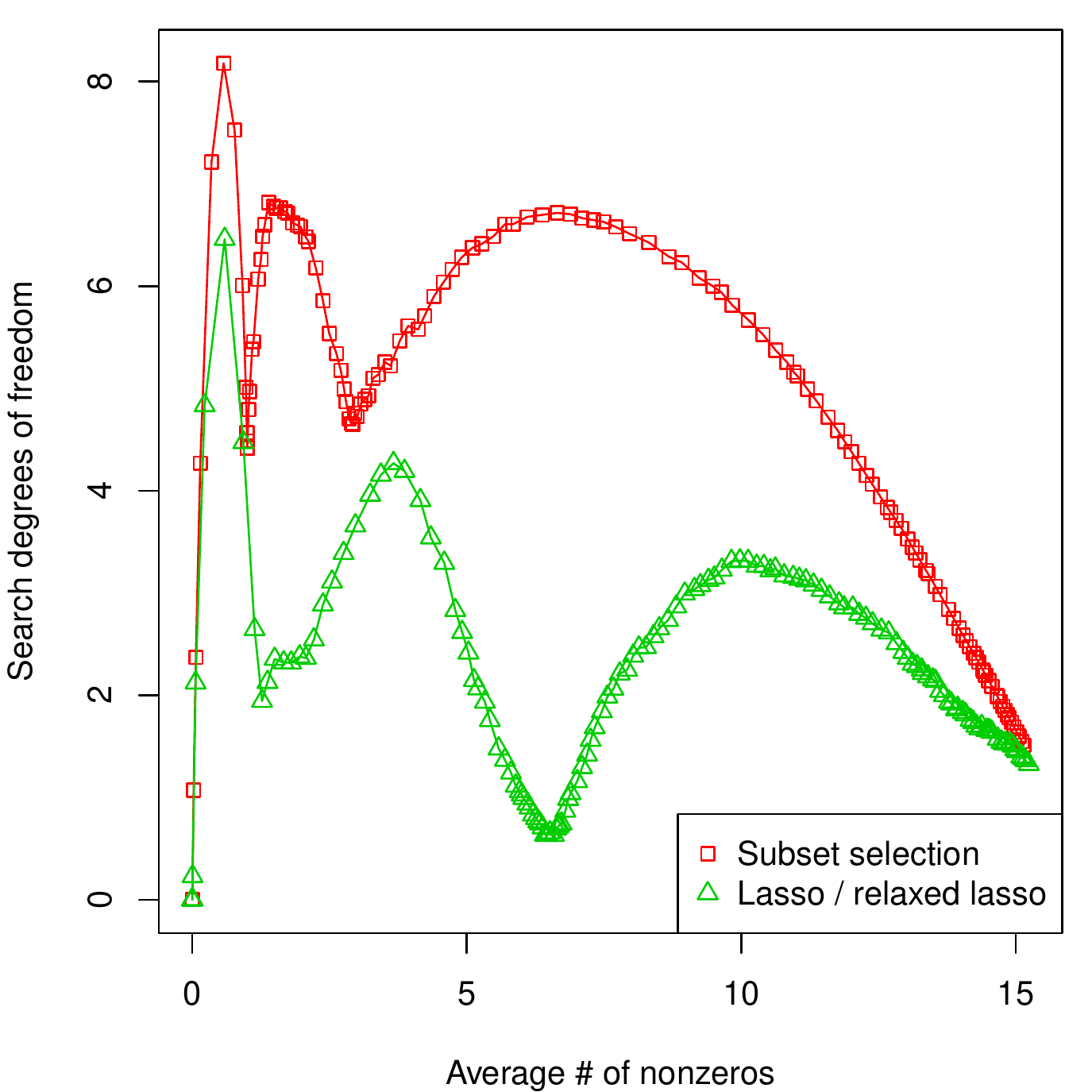} \\
\hline
Dense &
\includegraphics[width=\d]{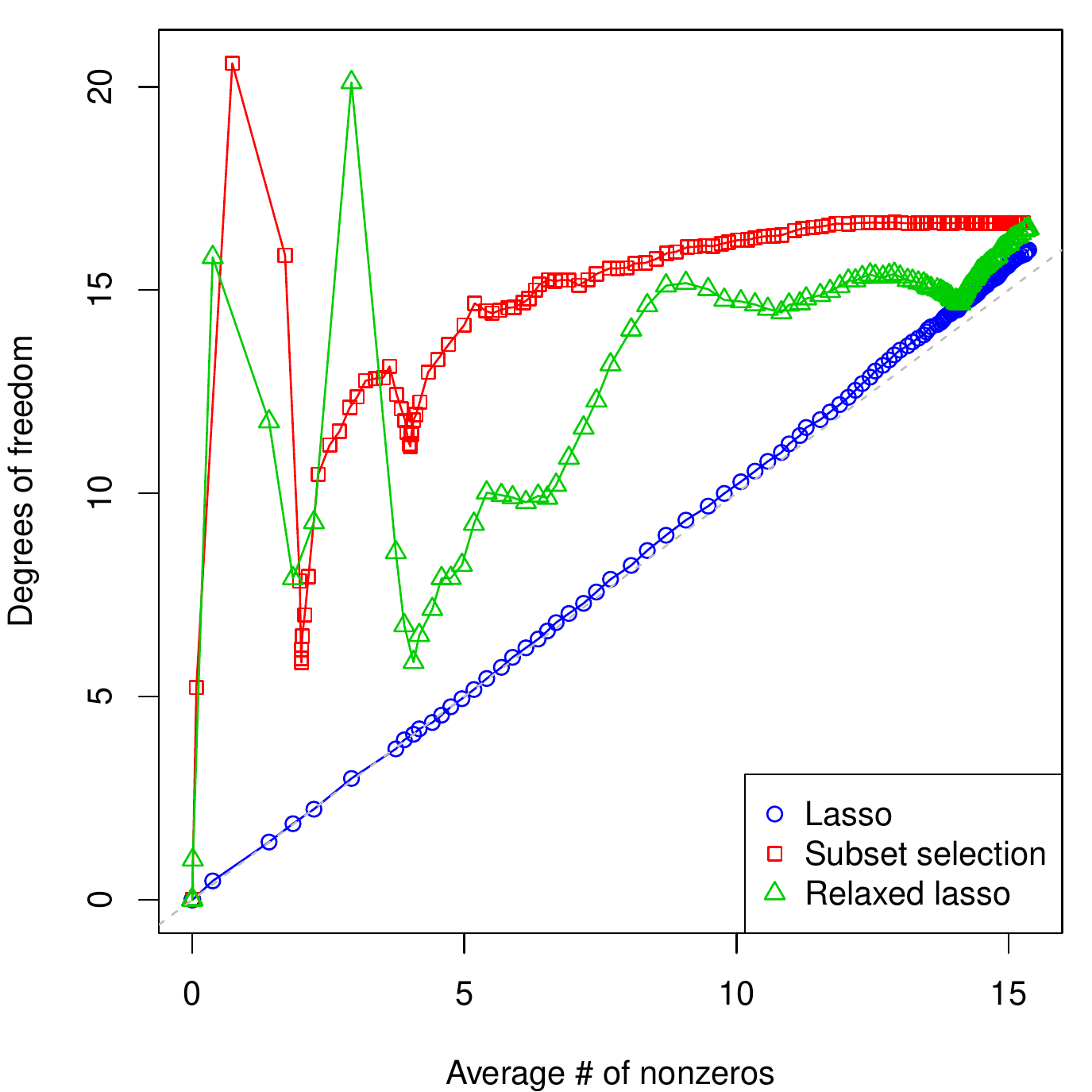} &
\includegraphics[width=\d]{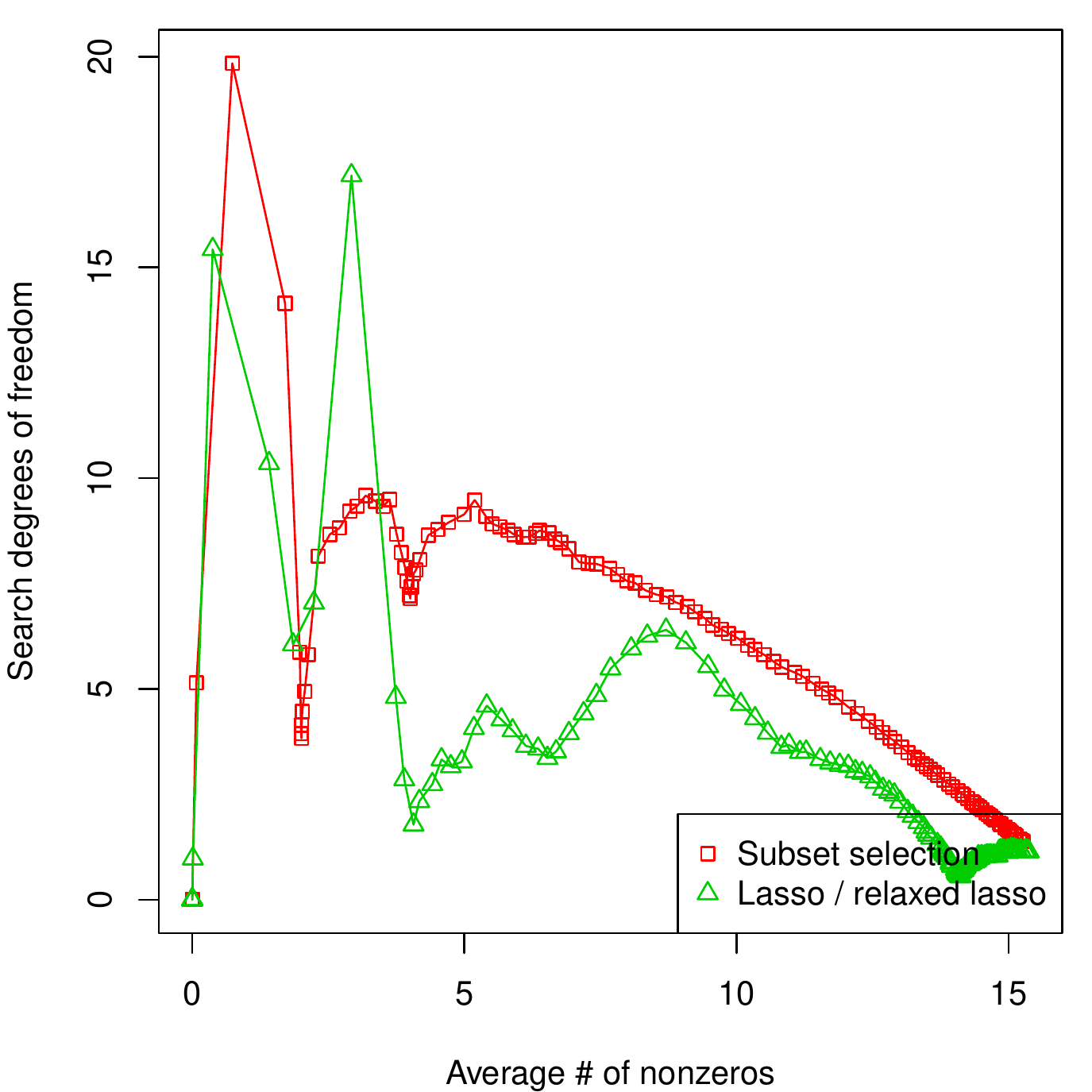} \\
\hline
\end{tabular}
\caption{\small\it A set of simulation results with $n=30$,
  $p=16$ (we are confined to such a small setup because of the
  exponential computational complexity of subset selection).  The
  rows of $X$ were drawn i.i.d.\ from $N(0,\Sigma)$, 
  where $\Sigma$ is block diagonal with two equal sized ($8\times 8$) 
  blocks $\Sigma_1, \Sigma_2$.  All diagonal entries of
  $\Sigma_1,\Sigma_2$ were set to 1, and the off-diagonal entries were
  drawn uniformly between 0.4 and 0.9.  We considered three cases for
  the true mean $\mu=X\beta^*$: null ($\beta^*=0$), sparse ($\beta^*$
  is supported on 3 variables in the first block and 1 in the second,
  with all nonzero components equal to 1), and dense ($\beta^*$ has
  all components equal to 1).
  In all cases, we drew $y$ around $\mu$ with independent standard
  normal noise, for a total of 100 repetitions.  Overall, the search
  degrees of freedom of subset selection appears to be larger than
  that of the lasso, but at times the latter can rival the former in
  magnitude, especially for small active sets, and in the dense signal 
  case.}   
\label{fig:lsr2}
\end{figure}

Unfortunately, this is not an easy task, since direct 
calculation of the relevant quantities---the degrees of freedom of
best subset selection 
and the relaxed lasso---is not tractable for a general $X$.  In
cases such as these, one usually turns to Stein's formula as an
alternative for calculating degrees of freedom; e.g., the result in
Theorem  \ref{thm:lassodf} is derived using Stein's formula.  But
Stein's formula only applies to continuous (and almost differentiable) 
fitting procedures $f=f(y)$, and neither the best subset selection nor
the relaxed lasso fit is continuous in $y$.   The next section,
therefore, is focused on extending Stein's result to discontinuous
functions.

\section{An extension of Stein's formula}
\label{sec:steinex}

This section considers Stein's formula \citep{stein}, and presents an
extension that yields an alternative derivation of the degrees of
freedom results in Section \ref{sec:orthx}, as well as (potential)
insights into the empirical results in Section \ref{sec:genx}.  In his
remarkable paper, Stein
studies the problem of estimating the mean of a multivariate normal
distribution, with a spherical covariance matrix, under the usual 
squared error loss. 
The main result is an unbiased estimate of the associated risk
for a large class of estimates of the mean. At the root of
Stein's arguments lies the following lemma.

\begin{lemma}[\citealt{stein}]
\label{lem:steinuni}
Let $Z \sim N(0,1)$. Let $f : \R \rightarrow
\R$ be absolutely continuous, with derivative $f'$. Assume that
$\E|f'(Z)| < \infty$. Then   
\begin{equation*}
\E[Zf(Z)] = \E[f'(Z)].
\end{equation*}
\end{lemma}

In its own right, this lemma (along with a converse statement, which
is also true) has a number of important applications
that span various areas of probability and statistics. 
For our purposes, the most relevant application is
an alternative and highly useful formula for computing degrees of
freedom.  This is given by extending the above lemma to
a setting in which the underlying normal distribution has an arbitrary
mean vector and variance, and is also multivariate. 

\begin{lemma}[\citealt{stein}]
\label{lem:steinmulti}
Let $X \sim N(\mu,\sigma^2 I)$, for some fixed $\mu \in \R^n$ and
$\sigma^2>0$. 
Let $g : \R^n \rightarrow \R$ be continuous and almost
differentiable, and write $\nabla g = (\partial g_1/\partial x_1,
\ldots \partial g_n/x_n)$ for the vector of partial
derivatives. Assume that $\E\|\nabla g(X)\|_2 < \infty$. Then 
\begin{equation}
\label{eq:steinmulti}
\frac{1}{\sigma^2} \E[(X-\mu)g(X)] = \E[\nabla g (X)].
\end{equation}
\end{lemma}

We will delay the definition of almost differentiability until a
little while later, but the eager reader can look ahead to
Definition \ref{def:ad}. Putting aside any concerns about regularity
conditions, the result in \eqref{eq:steinmulti} looks like a statement
about degrees of freedom.   
To complete the connection, consider a function $f : \R^n \rightarrow 
\R^n$ giving the fit $\hmu = f(y)$, and assume the usual normal 
model $y \sim N(\mu,\sigma^2 I)$. Let $f_i : \R^n \rightarrow \R$ be
the $i$th coordinate function of $f$. If $f_i$ satisfies the appropriate
conditions (continuity and almost differentiability), then we can
apply Lemma \ref{lem:steinmulti} with $X=y$ and $g=f_i$, take the $i$th
equality in \eqref{eq:steinmulti}, and sum over $i$ to give  
\begin{equation}
\label{eq:steindf}
\df(f) = \frac{1}{\sigma^2} \sum_{i=1}^n \Cov\big(f_i(y),y_i \big) = 
\E\bigg[ \sum_{i=1}^n \frac{\partial f_i}{\partial y_i} (y)\bigg],
\end{equation} 
where $\partial f_i/\partial y_i$ denotes the partial derivative of
$f_i$ with respect to its $i$th variable. This is known as  
{\it Stein's formula for degrees of freedom}. It can be very useful, 
because in some cases the divergence $\sum_{i=1}^n \partial
f_i/\partial y_i$ on the right-hand side of \eqref{eq:steindf} can be
computed explicitly, which yields an unbiased estimate of degrees
of freedom.  This is true, e.g., of the lasso fit, and as a conrete
illustration, we prove the result in Theorem \ref{thm:lassodf}, using
Stein's formula, in the appendix.

Useful as it can be, Stein's formula \eqref{eq:steindf} is not
universally applicable. There are several ways to 
break its assumptions; our particular interest is in fitting
procedures that are discontinuous in $y$. For example, we showed in
the proof of Theorem \ref{thm:subsetdf} that, when $X$ is orthogonal,
the subset selection solution is given by hard thresholding
$X^T y$ at the level $t=\sqrt{2\lambda}$. The hard thresholding
function $H_t$ is clearly discontinuous: each one of its coordinate
functions is discontinuous at $t$ and $-t$. We therefore derive a
modest extension of   
Stein's formula that allows us to deal with a certain class of
(well-behaved) discontinuous functions. We begin with the univariate
case, and then move on to the multivariate case. 

\subsection{An extension of Stein's univariate lemma}
\label{sec:steinexuni}

We consider functions $f: \R \rightarrow \R$ that are absolutely 
continuous on a partition of $\R$. Formally:

\begin{definition}
\label{def:pac}
We say that a function $f : \R \rightarrow \R$ is 
{\it piecewise absolutely continuous}, or 
{\it p-absolutely continuous}, if there exist points  
$\delta_1 < \delta_2 < \ldots < \delta_m$
such that $f$ is absolutely continuous on each one of the open
intervals   
$(-\infty, \delta_1), (\delta_1,\delta_2), \ldots (\delta_m,\infty)$. 
\end{definition}

For a p-absolutely continuous function $f$, we write
$\cD(f)=\{\delta_1,\ldots \delta_m\}$ for its discontinuity
set. Furthermore, note that such a function $f$ has a derivative $f'$
almost everywhere [because it has a derivative almost everywhere on
each of the intervals $(-\infty, \delta_1), (\delta_1,\delta_2), \ldots
(\delta_m,\infty)$]. We will simply refer to $f'$ as its derivative. 
Finally, we use the following helpful notation for one-sided limits, 
\begin{equation*}
f(x)_+ = \lim_{t \downarrow x} f(t) 
\;\;\; \text{and} \;\;\;
f(x)_- = \lim_{t \uparrow x} f(t).
\end{equation*}

We now have the following extension of Stein's
univariate lemma, Lemma \ref{lem:steinuni}. 

\begin{lemma}
\label{lem:steinexuni}
Let $Z \sim N(0,1)$. Let $f : \R \rightarrow
\R$ be p-absolutely continuous, 
and have a discontinuity set
$\cD(f) = \{\delta_1,\ldots \delta_m\}$.  Let $f'$ be its derivative, 
and assume that $\E|f'(Z)|<\infty$. Then  
\begin{equation*}
\E[Z f(Z)] = \E[f'(Z)] + \sum_{k=1}^m \phi(\delta_k) 
\big[f(\delta_k)_+-f(\delta_k)_-\big].
\end{equation*}
\end{lemma}

The proof is similar to Stein's proof of Lemma \ref{lem:steinuni}, and 
is left to the appendix, for readability.  It is straightforward to
extend this result to a nonstandard normal distribution. 

\begin{corollary}
\label{cor:steinexuni}
Let $X \sim N(\mu,\sigma^2)$.
Let $h : \R \rightarrow \R$ be p-absolutely continuous,
with discontinuity set $\cD(h) = \{\delta_1,\ldots \delta_m\}$,
and derivative $h'$ satisfying $\E|h'(X)|<\infty$. Then 
\begin{equation*}
\frac{1}{\sigma^2} \E[(X-\mu)h(X)] = \E[h'(X)]
+ \frac{1}{\sigma} \sum_{k=1}^m
\phi\bigg(\frac{\delta_k-\mu}{\sigma}\bigg)
\big[h(\delta_k)_+ - h(\delta_k)_-\big].
\end{equation*}
\end{corollary}

With this extension, we can immediately say something about degrees of   
freedom, though only in a somewhat restricted setting. Suppose that  
$f : \R^n \rightarrow \R^n$ provides the fit $\hmu=f(y)$, and that
$f$ is actually univariate in each coordinate, 
\begin{equation*}
f(y) = \big( f_1(y_1),\ldots f_n(y_n) \big).
\end{equation*}
Suppose also that each coordinate function $f_i:\R \rightarrow \R$ is
p-absolutely continuous. We can apply Corollary
\ref{cor:steinexuni} with $X=y_i$ and $h=f_i$, and sum over $i$ to
give 
\begin{align}
\nonumber
\df(f) &= \frac{1}{\sigma^2}\sum_{i=1}^n \Cov\big(f_i(y_i),y_i\big) \\ 
\label{eq:steinexuni}
&= \sum_{i=1}^n \E[f_i'(y_i)] + \frac{1}{\sigma} \sum_{i=1}^n 
\sum_{\delta \in \cD(f_i)} \phi\bigg(\frac{\delta-\mu_i}{\sigma}\bigg)
\big[f_i(\delta)_+ - f_i(\delta)_-\big].
\end{align}
The above expression provides an alternative way of proving the
result on the degrees of freedom of hard thresholding, which was given
in Lemma \ref{lem:htdf}, the critical lemma for deriving the degrees of
freedom of both best subset selection and the relaxed lasso for
orthogonal predictors, Theorems 
\ref{thm:subsetdf} and \ref{thm:relaxdf}. We step through this proof
next. 

\begin{proof}[Alternate proof of Lemma \ref{lem:htdf}]
For $f(y)=H_t(y)$, the $i$th coordinate function is 
\begin{equation*}
f_i(y_i) = [H_t(y_i)]_i = y_i \cdot 1\{|y_i| \geq t\}, 
\end{equation*}
which has a discontinuity set $\cD(f_i) = \{-t,t\}$. The second term in 
\eqref{eq:steinexuni} is hence 
\begin{equation*}
\frac{1}{\sigma} \sum_{i=1}^n 
\bigg[
\phi\bigg(\frac{t-\mu_i}{\sigma}\bigg)\cdot(t - 0) +
\phi\bigg(\frac{-t-\mu_i}{\sigma}\bigg)\cdot(0 - -t) 
\bigg] = \frac{t}{\sigma}
\sum_{i=1}^n \bigg[\phi\bigg(\frac{t-\mu_i}{\sigma}\bigg) +  
\phi\bigg(\frac{t+\mu_i}{\sigma}\bigg)\bigg],
\end{equation*}
while the first term is simply
\begin{equation*}
\sum_{i=1}^n \E[1\{|y_i| \geq t\}] = \E|\cA_t|.
\end{equation*}
Adding these together gives
\begin{equation*}
\df(H_t) = \E|\cA_t| + \frac{t}{\sigma}
\sum_{i=1}^n \bigg[\phi\bigg(\frac{t-\mu_i}{\sigma}\bigg) +  
\phi\bigg(\frac{t+\mu_i}{\sigma}\bigg)\bigg],
\end{equation*}
precisely the conclusion of Lemma \ref{lem:htdf}. 
\end{proof}

\subsection{An extension of Stein's multivariate lemma}
\label{sec:steinexmulti}

The degrees of freedom result \eqref{eq:steinexuni} applies to
functions $f$ for which the $i$th component function $f_i$
depends only on the $i$th component of the input, $f_i(y) = f_i(y_i)$,
for $i=1,\ldots n$.   Using this result,
we could compute the degrees of freedom of the best subset  
selection and relaxed lasso fits in the orthogonal predictor matrix
case.  Generally speaking, however, we cannot use this result outside
of the orthogonal setting, due to the requirement on $f$ that
$f_i(y)=f_i(y_i)$, $i=1,\ldots n$.  Therefore, in the hope of
understanding degrees of freedom for procedures like best subset
selection and the relaxed lasso in a broader context, we derive an
extension of Stein's multivariate lemma.     

Stein's multivariate lemma, Lemma \ref{lem:steinmulti}, is concerned
with functions $g : \R^n \rightarrow \R$ that are continuous and
almost differentiable.   
Loosely speaking, the concept of almost differentiability is really a
statement about absolute continuity. In words, a function is said to
be almost differentiable if it is absolutely continuous on almost
every line parallel to the coordinate axis (this notion is different,
but equivalent, to that given by Stein). 
Before translating this mathematically, we introduce some
notation. Let us write $x = (x_i,x_{-i})$ to emphasize that $x \in
\R^n$ is determined by its 
$i$th component $x_i \in \R$ and the other $n-1$ components 
$x_{-i} \in \R^{n-1}$. For $g : \R^n \rightarrow \R$, we let
$g(\,\cdot\,,x_{-i})$ denote $g$ as a function of the $i$th component
alone, with all others components fixed at the value $x_{-i}$. We can
now formally define almost differentiability:

\begin{definition}
\label{def:ad}
We say that a function $g:\R^n \rightarrow \R$ is 
{\it almost differentiable} if for every $i=1,\ldots n$ and Lebesgue
almost every $x_{-i} \in \R^{n-1}$, the function 
$g(\,\cdot\,,x_{-i}) : \R \rightarrow \R$ is absolutely continuous.  
\end{definition}

Similar to the univariate case, we
propose a relaxed continuity condition. Namely:

\begin{definition}
\label{def:pad}
We say that a function $g:\R^n \rightarrow \R$ is
{\it p-almost differentiable} if for every $i=1,\ldots n$ and Lebesgue 
almost every $x_{-i} \in \R^{n-1}$, the function 
$g(\,\cdot\,,x_{-i}) : \R \rightarrow \R$ is p-absolutely continuous.  
\end{definition}

Note that a function $g$ that is p-almost differentiable has
partial derivatives almost everywhere, and we write the collection as 
$\nabla g = 
(\partial g/\partial x_1,\ldots \partial g/\partial x_n)$.\footnote{Of
  course, this does not necessarily mean that $g$ has a well-defined
  gradient, and so, cumbersome as it may read, we are careful about
  referring to $\nabla g$ as the vector of partial derivatives,
  instead of the gradient.}
Also, when dealing with $g(\,\cdot\,,x_{-i})$, the function $g$
restricted to its $i$th variable with all others fixed at $x_{-i}$, we
write its one-sided limits as
\begin{equation*}
g(x_i,x_{-i})_+ = \lim_{t \downarrow x_i} g(t,x_{-i}) \;\;\; \text{and} \;\;\;
g(x_i,x_{-i})_- = \lim_{t \uparrow x_i} g(t,x_{-i}).
\end{equation*}
We are now ready to present our extension of Stein's multivariate
lemma. 

\begin{lemma}
\label{lem:steinexmulti}
Let $X \sim N(\mu,\sigma^2 I)$, for some fixed $\mu \in \R^n$
and $\sigma^2>0$.  Let $g : \R^n \rightarrow \R$ be p-almost
differentiable, with vector of partial derivatives $\nabla g =
(\partial g/\partial x_1,\ldots \partial g/\partial x_n)$. Then,
for $i=1,\ldots n$, 
\begin{equation*}
\frac{1}{\sigma^2} \E[(X_i-\mu_i) g(X)] = \E\bigg[\frac
{\partial g}{\partial x_i}(X)\bigg] +  
\frac{1}{\sigma} 
\E\Bigg[\sum_{\delta \in \cD(g(\,\cdot\,,X_{-i}))}
\phi\bigg(\frac{\delta-\mu_i}{\sigma}\bigg)
\big[g(\delta,X_{-i})_+ - g(\delta,X_{-i})_-\big]\Bigg], 
\end{equation*}
provided that $\E|\partial g/\partial x_i\,(X)| < \infty$ and 
\begin{equation*}
\E\,\Bigg|\sum_{\delta \in \cD(g(\,\cdot\,,X_{-i}))}
\phi\bigg(\frac{\delta-\mu_i}{\sigma}\bigg)
\big[g(\delta,X_{-i})_+ - g(\delta,X_{-i})_-\big]\Bigg|
< \infty.
\end{equation*}
\end{lemma}

Refer to the appendix for the proof, which utilizes a conditioning
argument to effectively reduce the multivariate setup to a
univariate one, and then invokes Lemma \ref{lem:steinexuni}.

The above lemma, Lemma \ref{lem:steinexmulti}, leads to our most
general extension of Stein's formula for degrees of freedom.
Let $f : \R^n \rightarrow \R^n$ be a fitting procedure, as in 
$\hmu = f(y)$, and let $f(y) = (f_1(y),\ldots
f_n(y))$. Consider Lemma \ref{lem:steinexmulti} applied to the $i$th
coordinate function, so that $X=y$ and $g=f_i$. Provided that each
$f_i$ is p-almost differentiable and satisfies the regularity
conditions 
\begin{equation}
\label{eq:reg}
\E\bigg|\frac{\partial f_i}{\partial y_i}(y)\bigg| < \infty  
\;\;\;\text{and}\;\;\;
\E\,\Bigg|\sum_{\delta \in \cD(f_i(\,\cdot\,,y_{-i}))}
\phi\bigg(\frac{\delta-\mu_i}{\sigma}\bigg)
\big[f_i(\delta,y_{-i})_+ - f_i(\delta,y_{-i})_-\big]\Bigg|
< \infty,
\end{equation}
we can take the $i$th equality in the lemma, and sum over $i$ to  
give 
\begin{align}
\nonumber
\df(f) &= \frac{1}{\sigma^2}\sum_{i=1}^n \Cov\big(f_i(y),y_i\big) \\
\label{eq:steinexmulti}
&= \sum_{i=1}^n \E\bigg[\frac{\partial f_i}{\partial y_i}(y)\bigg]  
+ \frac{1}{\sigma} \sum_{i=1}^n 
\E\Bigg[\sum_{\delta \in \cD(f_i(\,\cdot\,,y_{-i}))}
\phi\bigg(\frac{\delta-\mu_i}{\sigma}\bigg)
\big[f_i(\delta,y_{-i})_+ - f_i(\delta,y_{-i})_-\big]\Bigg].
\end{align}
Even if we assume that \eqref{eq:steinexmulti} is applicable to subset  
selection and the relaxed lasso with arbitrary predictors $X$, the
discontinuity sets---and hence the second term in
\eqref{eq:steinexmulti}---seem to be quite difficult to calculate in these
cases. In other words, unfortunately, the formula \eqref{eq:steinexmulti}
does not seem to provide an avenue for exact computation of the
degrees of freedom of subset selection or the relaxed lasso in
general.  However, it may still help us understand these
quantities, as we discuss next.

\subsection{Potential insights from the multivariate Stein extension} 


For both of the best subset selection and relaxed lasso fitting
procedures, one can show that the requisite regularity conditions
\eqref{eq:reg} indeed hold, which makes the extended Stein formula 
\eqref{eq:steinexmulti} valid.  Here we briefly outline a  
geometric interpretation for these fits, and describe how it can be
used to understand their discontinuity sets, and the formula in
\eqref{eq:steinexmulti}.  For an argument of a similar kind (and one   
given in more rigorous detail), see \citet{lassodf2}.  

In both cases, we can decompose $\R^n$ into a finite union of disjoint
sets, $\R^n = \cup_i^m U_i$, with each $U_i$ being polyhedral for the 
relaxed lasso, and each $U_i$ an intersection of quadratic sublevel
sets for subset selection. 
The relaxed lasso and best subset selection now share the property 
that, on the relative interior of each set $U_i$ in their respective
decompositions, the fit is just a linear projection map, and assuming
that $X$ has columns in general position, this is just the projection
map onto the column space of $X_\cA$ for some fixed set $\cA$.
Hence the discontinuity set of the fitting procedure in each case is  
contained in $\cup_{i=1}^m \mathrm{relbd}(U_i)$, which has measure
zero.  In other words, the active set is locally constant for almost
every $y\in\R^n$, and only for $y$ crossing the relative
boundary of some $U_i$ does it change.  From this, we
can verify the appropriate regularity conditions in \eqref{eq:reg}.

As $f$ for the relaxed lasso and subset selection is the locally
linear projection map $f(y)=P_\cA y$, almost everywhere in $y$, the
first term $\sum_{i=1}^n \E[\partial f_i(y)/\partial y_i]$ in
\eqref{eq:steinexmulti} is simply $\E|\cA|$.  The second term, then, 
exactly coincides with the search degrees of freedom of
these procedures. (Recall that the same breakdown occurred when using  
the univariate Stein extension to derive the degrees of freedom of
hard thresholding, in Section \ref{sec:steinexuni}.)  This suggests
a couple potential insights into degrees of freedom and search
degrees of freedom that may be gleaned from the extended Stein formula
\eqref{eq:steinexmulti}, which we discuss below.  

\begin{itemize}
\item {\it Positivity of search degrees of freedom.}  If one could
  show that
\begin{equation}
\label{eq:discsign}
f_i(\delta,y_{-i})_+ - f_i(\delta,y_{-i})_- > 0
\end{equation}
for each discontinuity point $\delta \in \cD(f_i(\cdot,y_{-i}))$,
almost every $y_{-i} \in \R^n$, and each $i=1,\ldots n$, then this
would imply that the second term in \eqref{eq:steinexmulti} is
positive.  For the relaxed lasso and subset selection fits,
this would mean that the search degrees of freedom term is always 
positive, i.e., 
the total degrees of freedom of these procedures is
always larger than the (expected) number of selected
variables.  In words, the condition in \eqref{eq:discsign} says that 
the $i$th fitted value, at a point of discontinuity, can only increase 
as the $i$th component of $y$ increases.  Note that this is a
sufficient but not necessary condition for positivity of search
degrees of freedom.

\item {\it Search degrees of freedom and discontinuities.}  The fact
  that the second term in \eqref{eq:steinexmulti} gives the search
  degrees of freedom of the best subset selection and the relaxed
  lasso fits tells us that the search degrees of freedom of a
  procedure is intimately related to its discontinuities over $y$.
  At a high level: the greater the number of discontinuities, the
  greater the magnitude of these discontinuities, and the closer
  they occur to the true mean $\mu$, the greater the search degrees of
  freedom.   

  This may provide some help in understanding 
  the apparent (empirical) differences in search
  degrees of freedom between the relaxed lasso and best subset
  selection fits under correlated setups, as seen in Section
  \ref{sec:genx}.  The particular discontinuities of concern in 
  \eqref{eq:steinexmulti} arise from fixing all but 
  $i$th component of the outcome at $y_{-i}$, and examining the $i$th
  fitted value $f_i(\cdot,y_{-i})$ as a function of its $i$th
  argument.  One might expect that this function $f_i(\cdot,y_{-i})$
  generally exhibits more points of discontinuity for best subset
  selection compared to the relaxed lasso, due to the
  more complicated boundaries of the elements $U_i$ in the
  active-set-determining decomposition described above (these
  boundaries are piecewise quadratic for best subset selection, and 
  piecewise linear for the relaxed lasso).  This is in line with the
  general trend of subset selection displaying a larger search
  degrees of freedom than the relaxed lasso.

  But, as demonstrated in Figure
  \ref{fig:lsr2}, something changes for large values of
  $\lambda$ (small active sets, on average), and for $\mu=X\beta^*$
  with a sparse or (especially) dense true coefficient vector
  $\beta^*$; we saw that the search degrees of freedom of both the
  relaxed lasso and best subset selection fits can grow very large in
  these cases.  
  Matching search degrees of freedom to the second term in
  \eqref{eq:steinexmulti}, therefore, we infer that both fits must
  experience major discontinuities here (and these are somehow
  comparable overall, when measured in number, magnitude, and
  proximity to $\mu$).  This makes sense, especially when we think of
  taking $\lambda$ large enough so that these procedures are forced to 
  select an active set that is strictly contained in the true support
  $\cA^*=\supp(\beta^*)$; different values of $y$, quite close to
  $\mu=X\beta^*$, will make different subsets of $\cA^*$ look more or
  less appealing according to the criterions in \eqref{eq:lasso},
  \eqref{eq:subset}.  
\end{itemize}

\subsection{Connection to Theorem 2 of \citet{dfnonlin}}

After completing this work, we discovered the independent and
concurrent work of \citet{dfnonlin}.  These authors propose an
interesting and completely different geometric approach to studying 
the degrees of freedom of a metric projection estimator
\begin{equation*}
f(y) \in \argmin_{u \in K} \, \|y-u\|_2^2,
\end{equation*}
where the set $K \subseteq \R^n$ can be nonconvex.  
Their Theorem 2 gives a decomposition for degrees of freedom
that possesses an intriguing tie to ours in \eqref{eq:steinexmulti}.  
Namely, these authors show that the degrees of freedom of any metric
projection estimator $f$ an be expressed as its expected divergence
plus an ``extra'' term, this term being the integral of the normal
density with respect to a singular measure (dependent on $f$). 
Equating this with our expression  
in \eqref{eq:steinexmulti}, we see that the two forms of ``extra''
terms must match---i.e., our second term in \eqref{eq:steinexmulti},
defined by a sum over the discontinuities of the projection
$f$, must be equal to their integral.   

This has an immediate implication for the projection operator onto the
$\ell_0$ ball of radius $k$, i.e., the best subset selection estimator
in constrained form: the search degrees of freedom here must be
nonnegative (as the integral of a density with respect to a measure is
always nonnegative).  The decomposition of \citet{dfnonlin} hence 
elegantly proves that the best subset selection fit, constrained to
have $k$ active variables, attains a degrees of freedom larger than
or equal $k$.   However, as far as we can tell, their Theorem 2 does
not apply to best subset selection in Lagrange form, the estimator
considered in our paper, since it is limited to metric
projection estimators.  To be clear, our extension of Stein's formula
in \eqref{eq:steinexmulti} is not restricted to any particular form of 
fitting procedure $f$ (though we do require the regularity conditions
in \eqref{eq:reg}).  

We find the connections between our work and theirs fascinating, and
hope to understand them more deeply in the future.

\section{Discussion}
\label{sec:discuss}

In this work, we explored the degrees of freedom of best
subset selection and the relaxed lasso (the procedure that performs
least squares on the active set returned by the lasso).  We 
derived exact expressions for the degrees of freedom of these fitting
procedures with orthogonal predictors $X$, and investigated by
simulation their degrees of freedom for correlated predictors.  We
introduced a new concept, search degrees of freedom, which intuitively 
measures the amount of degrees of freedom expended by an adaptive
regression procedure in merely constructing an active set of variables
(i.e., not counting the degrees of freedom attributed to estimating
the active coefficients).  Search degrees of freedom has a precise
definition for any regression procedure.  For subset selection
and the relaxed lasso, this reduces to the (total) degrees of
freedom minus the expected number of active variables; for the lasso,
we simply equate its search degrees of freedom with that of the
relaxed lasso, since these two procedures have the exact same search
step.   

The last section of this paper derived an extension of Stein's formula
for discontinuous functions.  This was motivated by the hope
that such a formula could provide an alternative lens from which we
could view degrees of freedom for discontinuous fitting procedures
like subset selection and the relaxed lasso.  The application of this
formula to these fitting procedures is not easy, and our grasp of 
the implications of this formula for degrees of freedom is only
preliminary. There is much work to be done, but we are
hopeful that our extension of Stein's result will prove useful for
understanding degrees of freedom and search degrees of freedom, and
potentially, for other purposes as well.

\section*{Acknowledgements}

The idea for this paper was inspired by a conversation with Jacob 
Bien.  We thank Rob Tibshirani for helpful feedback
and encouragement, and the editors and
referees who read this paper and gave many useful comments
and references.

\appendix
\section{Proofs}

\subsection{Proof of Lemma \ref{lem:htdf}}

By definition,
\begin{align}
\nonumber
\df(H_t) &= \frac{1}{\sigma^2} 
\sum_{i=1}^n \Cov \big([H_t(y)]_i,y_i\big) \\ 
\nonumber
&= \frac{1}{\sigma^2} \sum_{i=1}^n \E\Big[
y_i(y_i-\mu_i) \Big(1\{y_i \geq t\} + 1\{y_i \leq -t\}\Big)\Big] \\  
\label{eq:dflast}
&= \frac{1}{\sigma^2} \sum_{i=1}^n \E\Big[
(z_i+\mu_i)z_i \Big(1\{z_i \geq t-\mu_i\} + 
1\{z_i \leq -t-\mu_i\}\Big)\Big],
\end{align}
where $z = y-\mu \sim N(0,\sigma^2 I)$.  To compute the above, we note
the identities (the last two can be checked using integration by parts): 
\begin{align}
\label{eq:tm1}
\E\big[z_i \cdot 1\{z_i \leq a\}\big] &= -\sigma \phi(a/\sigma), \\
\label{eq:tm2}
\E\big[z_i \cdot 1\{z_i \geq b\}\big] &= \sigma \phi(b/\sigma), \\
\label{eq:tm3}
\E\big[z_i^2 \cdot 1\{z_i \leq a\}\big] &= -\sigma a \phi(a/\sigma) + 
\sigma^2 \Phi(a/\sigma), \\
\label{eq:tm4}
\E\big[z_i^2 \cdot 1\{z_i \geq b\}\big] &= \sigma b \phi(b/\sigma) +
\sigma^2 \big[1-\Phi(b/\sigma)\big],
\end{align}
where $\Phi$ denotes the standard normal cdf.
Plugging these in, the expression in \eqref{eq:dflast} becomes 
\begin{equation*}
\sum_{i=1}^n\bigg[
1-\Phi\bigg(\frac{t-\mu_i}{\sigma}\bigg)
+\Phi\bigg(\frac{-t-\mu_i}{\sigma}\bigg)
\bigg] + 
 \frac{t}{\sigma} \sum_{i=1}^n 
\bigg[\phi\bigg(\frac{t-\mu_i}{\sigma}\bigg) + 
\phi\bigg(\frac{-t-\mu_i}{\sigma}\bigg)\bigg],
\end{equation*}
and the first sum above is exactly
\begin{equation*}
\E\bigg[\sum_{i=1}^n \Big(1\{z_i \geq t-\mu_i\} +
1\{z_i \leq -t-\mu_i\}\Big)\bigg] = 
\E|\cA_t|,
\end{equation*} 
as desired.
\qedsymbol

\subsection{Proof of Theorem \ref{thm:subsetdf}}

As $X$ is orthogonal, the criterion in \eqref{eq:subset} can be
written as 
\begin{equation*}
\|y-X\beta\|_2^2 = \|X^T y - \beta\|_2^2 + c,
\end{equation*}
where $c$ is a constant, meaning that it does not depend on $\beta$. 
Hence we can rewrite the optimization problem in \eqref{eq:subset} as 
\begin{equation*}
\hbetasub \in \argmin_{\beta \in \R^p} \;  
\half\|X\T y-\beta\|_2^2 + \lambda\|\beta\|_0,
\end{equation*}
and from this it is not hard to see that the solution is
\begin{equation*}
\hbetasub = H_{\sqrt{2\lambda}}(X\T y),
\end{equation*}
hard thresholding the quantity $X\T y$ by the amount
$t=\sqrt{2\lambda}$. Finally, 
we note that  
\begin{equation*}
\df(X\hbetasub) = \tr\Big(\Cov(X\hbetasub,y)\Big) =
\tr\Big(\Cov(\hbetasub,X\T y)\Big),
\end{equation*}
because the trace operator is invariant under commutation of
matrices, and $X\T y \sim N(X\T \mu, \sigma^2 I)$. 
Applying Lemma \ref{lem:htdf} completes the proof. 
\qedsymbol

\subsection{Proof of Theorem \ref{thm:lassodf}}

We use several facts about the lasso without proof. These are
derived in, e.g., \citet{lassodf2} and \citet{lassounique}. 
For fixed $X,\lambda$, the lasso fit $f(y)=\hmulas(y)$ is continuous
and almost differentiable in each coordinate, so we can apply Stein's  
formula \eqref{eq:steindf}.  As $X$ has columns in general position,
there is a unique lasso solution \smash{$\hbetalas$}, and letting
$\cA$ denote its active set, and $s$ denote the signs of active
lasso coefficients,
\begin{equation*}
\cA=\supp(\hbetalas) 
\;\;\; \text{and} \;\;\;
s=\sign(\hbetalas_\cA), 
\end{equation*}
the fit can be expressed as
\begin{equation*}
\hmulas = X_\cA(X_\cA^T X_\cA)^{-1} X_\cA^T y - 
X_\cA (X_\cA^T X_\cA)^{-1} \lambda s.
\end{equation*}  
For almost every $y\in \R^n$, the set $\cA$ and vector $s$
are locally constant (with respect to $y$), and so they have zero 
derivative (with respect to $y$). Hence, for almost every $y$,
\begin{equation*}
\sum_{i=1}^n \frac{\partial \hmulas_i}{\partial y_i}(y) = 
\tr\Big(X_\cA (X_\cA^T X_\cA)^{-1} X_\cA\T \Big) =
|\cA|,
\end{equation*}
and taking an expectation gives the result. 
\qedsymbol

\subsection{Proof of Lemma \ref{lem:steinexuni}}

The result can be shown using integration by parts. We prove
it in a different way, mimicking Stein's proof of Lemma
\ref{lem:steinuni}, which makes the proof for the
multivariate case (Lemma \ref{lem:steinexmulti}) easier. We have 
\begin{align}
\nonumber
\E[f'(Z)] &= \int_{-\infty}^\infty f'(z) \phi(z)\,dz \\
\nonumber
&= \int_0^\infty f'(z) 
\bigg\{\int_z^\infty t\phi(t)\,dt\bigg\}\,dz -
\int_{-\infty}^0 f'(z) 
\bigg\{\int_{-\infty}^z t\phi(t)\,dt\bigg\}\,dz \\
\label{eq:fubini}
&= \int_0^\infty t\phi(t)
\bigg\{\int_0^t f'(z)\,dz\bigg\}\,dt -
\int_{-\infty}^0 t\phi(t)
\bigg\{\int_t^0 f'(z)\,dz\bigg\}\,dt.
\end{align}
The second equality follows from $\phi'(t)=-t\phi(t)$, and the
third is by Fubini's theorem. Consider the first term in
\eqref{eq:fubini}; as $f$ is absolutely continuous on each of the
intervals 
$(-\infty,\delta_1),(\delta_1,\delta_2), \ldots (\delta_m,\infty)$, 
the fundamental theorem of (Lebesgue) integral calculus gives 
\begin{equation*}
\int_0^t f'(z)\,dz = f(t)-f(0)-\sum_{k=1}^m
\big[f(\delta_k)_+ - f(\delta_k)_-\big] \cdot 1(0\leq\delta_k\leq t). 
\end{equation*}
Therefore
\begin{equation*}
\int_0^\infty t\phi(t)
\bigg\{\int_0^t f'(z)\,dz\bigg\}\,dt =
\int_0^\infty t\phi(t)\big[f(t)-f(0)\big]\,dt  
- \sum_{\delta_k\geq 0}
\big[f(\delta_k)_+ - f(\delta_k)_-\big]
\int_{\delta_k}^\infty t\phi(t) \,dt.
\end{equation*}
The second term in \eqref{eq:fubini} is similar, and putting these
together we get
\begin{multline*}
\E[f'(Z)] = \E[Zf(Z)]- \E[Z]f(0) -\sum_{\delta_k \geq 0}
\big[f(\delta_k)_+ - f(\delta_k)_-\big] \cdot 
\E\big[Z \cdot 1\{Z\geq \delta_k\}\big] \;+ \\  
\sum_{\delta_k < 0} \big[f(\delta_k)_+ - f(\delta_k)_-\big] 
\cdot \E\big[Z \cdot 1\{Z \leq \delta_k\}\big].
\end{multline*}
The result follows by noting that $\E[Z]=0$ and recalling the
identities \eqref{eq:tm1} and \eqref{eq:tm2}. 
\qedsymbol

\subsection{Proof of Corollary \ref{cor:steinexuni}}

Define $Z=(X-\mu)/\sigma$ and  
$f(z)=h(\sigma z + \mu)$, and apply Lemma \ref{lem:steinexuni}.
\qedsymbol

\subsection{Proof of Lemma \ref{lem:steinexmulti}}

We assume that $X \sim N(0,I)$, and then a similar standardization 
argument to that given in the proof of Corollary \ref{cor:steinexuni} 
can be applied here to prove the result for 
$X \sim N(\mu,\sigma^2 I)$.   
 
For fixed $X_{-i}$, the function $g(\,\cdot\,,X_{-i})$ is univariate. 
Hence, following the proof of Lemma \ref{lem:steinexuni}, and using 
the independence of $X_i$ and $X_{-i}$,
\begin{align*}
\E\bigg[\frac{\partial g}{\partial x_i}(X) \,\Big|\, X_{-i}\bigg]  
&= \int_{-\infty}^\infty \frac{\partial g}{\partial x_i} 
(z,X_{-i}) \phi(z)\,dz \\
&= \int_0^\infty \frac{\partial g}{\partial x_i} (z,X_{-i}) 
\bigg\{ \int_z^\infty t \phi(t) \, dt\bigg\} \, dz - 
\int_{-\infty}^0 \frac{\partial g}{\partial x_i} (z,X_{-i}) 
\bigg\{ \int_{-\infty}^z t \phi(t) \, dt\bigg\} \, dz \\
&= \int_0^\infty t\phi(t) \bigg\{ \int_0^t 
\frac{\partial g}{\partial x_i} (z,X_{-i})\,dz \bigg\} \, dt - 
\int_{-\infty}^0 t\phi(t) \bigg\{ \int_t^0 
\frac{\partial g}{\partial x_i} (z,X_{-i})\,dz \bigg\} \, dt.
\end{align*}
Consider the first term above. 
For almost every $X_{-i}$, the function $g(\,\cdot\,,X_{-i})$ is  
p-absolutely continuous, so the inner integral is 
\begin{equation*}
\int_0^t \frac{\partial g}{\partial x_i} (z,X_{-i}) \, dz =
g(t,X_{-i})-g(0,X_{-i}) -
\sum_{\delta \in \cD_i}
\big[g(z,X_{-i})_+ - g(z,X_{-i})_-\big] 
\cdot 1(0 \leq \delta \leq t),
\end{equation*}
where we have abbreviated $\cD_i=\cD(g(\,\cdot\,,X_{-i}))$. 
The next steps follow the corresponding arguments in the proof of
Lemma \ref{lem:steinexuni}, yielding
\begin{equation*}
\E\bigg[\frac{\partial g}{\partial x_i}(X)\,\Big|\,X_{-i}\bigg] =  
\E[X_i g(X)\,|\,X_{-i}] - \sum_{\delta \in \cD_i}  
\big[ g(\delta,X_{-i})_+ - g(\delta,X_{-i})_- \big]
\end{equation*} 
for almost every $X_{-i}$. Taking an expectation over $X_{-i}$
gives the result. 
\qedsymbol

\bibliographystyle{agsm}  
\bibliography{ryantibs}

\end{document}